\documentclass[10pt]{amsart}
\usepackage{eucal, mathrsfs}
\usepackage{amsmath,amssymb,amsthm,amsfonts,amscd,amsopn}
\usepackage{hyperref}
\usepackage{times}
\usepackage{xspace}
\usepackage{epsfig,epic,eepic,latexsym,color}
\usepackage[all]{xy}
\usepackage{comment, alltt} 
\usepackage{enumitem}
\xyoption{all}

\newcommand{\field}[1]{\mathbb #1} 
\newcommand{\mf}[1]{\mathfrak #1}
\newcommand{\mc}[1]{\mathcal #1}
\newcommand{\ms}[1]{\mathscr #1}
\newcommand{\widebar}[1]{\overline{#1}}

\newcommand{\LL}{\field L}

\newcommand{\Z}{\field Z}
\newcommand{\Q}{\field Q}

\newcommand{\WWW}{{\mc W}}
\newcommand{\WW}{{\sf W}}
\newcommand{\ww}{{\sf w}}
\newcommand{\setW}{{\sW}}
\newcommand{\schW}{{\sf W}}

\newcommand{\simto}{\stackrel{\sim}{\to}}
\newcommand{\eps}{\varepsilon}
\renewcommand{\phi}{\varphi}

\newcommand{\red}{\operatorname{red}}

\renewcommand{\hom}{\operatorname{Hom}}
\newcommand{\uhom}{\operatorname{\underline{Hom}}}
\newcommand{\shom}{\ms H\!om}

\newcommand{\spec}{\operatorname{Spec}}

\newcommand{\rspec}{\operatorname{\bf Spec}}

\renewcommand{\P}{\field P}
\newcommand{\A}{\field A}

\DeclareMathOperator{\Pic}{Pic}

\DeclareMathOperator{\pr}{pr}

\newcommand{\m}{\boldsymbol{\mu}}

\newcommand{\G}{\field G} 

\newcommand{\PGL}{\operatorname{PGL}}

\DeclareMathOperator{\coh}{\operatorname{Coh}}

\newcommand{\ch}{\operatorname{char}}

\DeclareMathOperator*{\tensor}{\otimes}

\newcommand{\inj}{\hookrightarrow}
\newcommand{\id}{\operatorname{id}}

\DeclareMathOperator{\aut}{\operatorname{Aut}}
\DeclareMathOperator{\isom}{\operatorname{Isom}}

\DeclareMathOperator{\V}{\bf V}

\DeclareMathOperator{\B}{\operatorname{B\!}}

\DeclareMathAlphabet{\smallchanc}{OT1}{pzc}%
                                 {m}{it}
\DeclareFontFamily{OT1}{pzc}{}
\DeclareFontShape{OT1}{pzc}{m}{it}%
             {<-> s * [1.100] pzcmi7t}{}
\DeclareMathAlphabet{\mathchanc}{OT1}{pzc}%
                                 {m}{it}

\newcommand{\mcL}{\mathchanc{L}}

\makeatletter
\setenumerate[1]{leftmargin=*,parsep=0em,itemsep=0.125em,topsep=0.125em}
\renewcommand\subsection{
  \renewcommand{\sfdefault}{pag}
  \@startsection{subsection}%
  {2}{0pt}{-\baselineskip}{.2\baselineskip}{\raggedright
    \sffamily\itshape\small
  }}
\renewcommand\section{
  \renewcommand{\sfdefault}{phv}
  \@startsection{section} %
  {1}{0pt}{\baselineskip}{.2\baselineskip}{\centering
    \sffamily
    \scshape
}}

\newtheoremstyle{bozont}{3pt}{3pt}%
     {\itshape}
     {}
     {\bfseries}
     {.}
     {.5em}
     {\thmname{#1}\thmnumber{ #2}\thmnote{ \rm #3}}
\newtheoremstyle{bozont-sf}{3pt}{3pt}%
     {\itshape}
     {}
     {\sffamily}
     {.}
     {.5em}
     {\thmname{#1}\thmnumber{ #2}\thmnote{ \rm #3}}
\newtheoremstyle{bozont-sc}{3pt}{3pt}%
     {\itshape}
     {}
     {\scshape}
     {.}
     {.5em}
     {\thmname{#1}\thmnumber{ #2}\thmnote{ \rm #3}}
\newtheoremstyle{bozont-remark}{3pt}{3pt}%
     {}
     {}
     {\scshape}
     {.}
     {.5em}
     {\thmname{#1}\thmnumber{ #2}\thmnote{ \rm #3}}
\newtheoremstyle{bozont-def}{3pt}{3pt}%
     {}
     {}
     {\bfseries}
     {.}
     {.5em}
     {\thmname{#1}\thmnumber{ #2}\thmnote{ \rm #3}}
\newtheoremstyle{bozont-reverse}{3pt}{3pt}%
     {\itshape}
     {}
     {\bfseries}
     {.}
     {.5em}
     {\thmnumber{#2.}\thmname{ #1}\thmnote{ \rm #3}}
\newtheoremstyle{bozont-reverse-sc}{3pt}{3pt}%
     {\itshape}
     {}
     {\scshape}
     {.}
     {.5em}
     {\thmnumber{#2.}\thmname{ #1}\thmnote{ \rm #3}}
\newtheoremstyle{bozont-reverse-sf}{3pt}{3pt}%
     {\itshape}
     {}
     {\sffamily}
     {.}
     {.5em}
     {\thmnumber{#2.}\thmname{ #1}\thmnote{ \rm #3}}
\newtheoremstyle{bozont-remark-reverse}{3pt}{3pt}%
     {}
     {}
     {\sc}
     {.}
     {.5em}
     {\thmnumber{#2.}\thmname{ #1}\thmnote{ \rm #3}}
\newtheoremstyle{bozont-def-reverse}{3pt}{3pt}%
     {}
     {}
     {\bfseries}
     {.}
     {.5em}
     {\thmnumber{#2.}\thmname{ #1}\thmnote{ \rm #3}}
\newtheoremstyle{bozont-def-newnum-reverse}{3pt}{3pt}%
     {}
     {}
     {\bfseries}
     {}
     {.5em}
     {\thmnumber{#2.}\thmname{ #1}\thmnote{ \rm #3}}
\theoremstyle{bozont}    
\newtheorem{proclaim}{Theorem}[section]
\newtheorem{thm}[proclaim]{Theorem}
\newtheorem{cor}[proclaim]{Corollary} 
\newtheorem{lem}[proclaim]{Lemma} 
\newtheorem{prop}[proclaim]{Proposition}

\theoremstyle{bozont-sc}
\newtheorem{proclaim-special}[proclaim]{\specialthmname}

\theoremstyle{bozont-remark}
\newtheorem{rem}[proclaim]{Remark}

\newtheorem*{SubHeading*}{\SubHeadingName}%
\newtheorem{SubHeading}[proclaim]{\SubHeadingName}
\newtheorem{sSubHeading}[equation]{\sSubHeadingName}
\newenvironment{demo}[1] {\def\SubHeadingName{#1}\begin{SubHeading}}
  {\end{SubHeading}}%
\newenvironment{demo-r}[1] {\def\SubHeadingName{#1}\begin{SubHeading-r}}
  {\end{SubHeading-r}}%
\newenvironment{subdemo-r}[1] {\def\sSubHeadingName{#1}  \begin{sSubHeading-r}}
  {\end{sSubHeading-r}} %
\newenvironment{demo*}[1] {\def\SubHeadingName{#1}\begin{SubHeading*}}
  {\end{SubHeading*}}%

\newtheorem{remark}[proclaim]{Remark}

\newtheorem{notn}[proclaim]{Notation}

\theoremstyle{bozont-def}    
\newtheorem{defn}[proclaim]{Definition}

\theoremstyle{bozont-reverse}

\theoremstyle{bozont-reverse-sc}
\newtheorem{proclaimr-special}[proclaim]{\specialthmname}
\newenvironment{proclaimspecialr}[1]%
{\def\specialthmname{#1}\begin{proclaimr-special}}%
{\end{proclaimr-special}}
\theoremstyle{bozont-remark-reverse}

\newtheorem{SubHeading-r}[proclaim]{\SubHeadingName}
\newtheorem{sSubHeading-r}[equation]{\sSubHeadingName}
\newtheorem{SubHeadingr}[proclaim]{\SubHeadingName}

\theoremstyle{bozont-def-newnum-reverse}    
\newtheorem{newnumr}[proclaim]{}
\theoremstyle{bozont-def-reverse}

\newtheorem{newnumspecial}[proclaim]{\specialnewnumname}

\numberwithin{equation}{proclaim}
\numberwithin{figure}{section} 

\newenvironment{enumerate-p}{
  \begin{enumerate}}
  {\setcounter{equation}{\value{enumi}}\end{enumerate}}

\newlength{\swidth}
\makeatother

\makeatletter

\newcommand\noin{\noindent}

\newcommand\dsize{\displaystyle}

\DeclareFontFamily{OMS}{rsfs}{\skewchar\font'60}
\DeclareFontShape{OMS}{rsfs}{m}{n}{<-5>rsfs5 <5-7>rsfs7 <7->rsfs10 }{}
\DeclareSymbolFont{rsfs}{OMS}{rsfs}{m}{n}
\DeclareSymbolFontAlphabet{\scr}{rsfs}

\newcommand{\sA}{\scr{A}}

\newcommand{\sF}{\scr{F}}

\newcommand{\sL}{\scr{L}}

\newcommand{\sN}{\scr{N}}
\newcommand{\sO}{\scr{O}}

\newcommand{\sT}{\scr{T}}

\newcommand{\sW}{\scr{W}}
\newcommand{\sX}{\scr{X}}
\newcommand{\sY}{\scr{Y}}
\newcommand{\sZ}{\scr{Z}}

\newcommand{\bN}{\mathbb{N}}

\newcommand{\bP}{\mathbb{P}}

\newcommand{\frC}{\mathfrak{C}}

\newcommand{\iso}{\simeq}    

\newcommand{\into}{\hookrightarrow}
\newcommand{\onto}{\twoheadrightarrow}

\newcommand{\col}{\colon}

\renewcommand{\P}{\mathbb{P}}

\DeclareMathOperator{\codim}{codim}
\DeclareMathOperator{\Hom}{Hom}

\DeclareMathOperator{\Isom}{Isom}

\DeclareMathOperator{\Proj}{{Proj}}

\DeclareMathOperator{\Sing}{{Sing}}
\DeclareMathOperator{\Spec}{{Spec}}

\def\ring#1.{\scr O_{#1}}
\def\map#1.#2.{#1 \to #2}

\def\longmap#1.#2.{#1 \longrightarrow #2}
\def\factor#1.#2.{\left. \raise 2pt\hbox{$#1$} \right/
\hskip -2pt\raise -2pt\hbox{$#2$}}
\def\pe#1.{\mathbb P(#1)}
\def\prs#1.{\mathbb P^{#1}}

\def\coh#1.#2.#3.{\dsize H^{#1}(#2,#3)}
\def\dimcoh#1.#2.#3.{\dsize h^{#1}(#2,#3)}
\def\ses#1.#2.#3.{0  \longrightarrow  #1   \longrightarrow 
 #2 \longrightarrow #3 \longrightarrow 0} 
\def\sesshort#1.#2.#3.{0
 \rightarrow #1 \rightarrow #2 \rightarrow #3 \rightarrow 0}
\def\shortses#1.#2.#3.{0  \rightarrow  #1   \rightarrow 
 #2  \rightarrow   #3 \rightarrow  0}
\def\Iff#1#2#3{
\hfil\hbox{\hsize =#1
\vtop{\noin #2}
\hskip.5cm 
\lower.5\baselineskip\hbox{$\Leftrightarrow$}\hskip.5cm
\vtop{\noin #3}}\hfil\medskip}
\def\myoplus#1.#2.{\underset #1 \to {\overset #2 \to \oplus}}

\newcommand\resto[1]{\vert_{#1}}
\newcommand\wb{\mathsf b}
\newcommand\rf{\kappa}
\newcommand\sncd{{simple} normal crossings divisor\xspace}

\makeatother

\title{Boundedness of families of canonically polarized manifolds: A
higher dimensional analogue of Shafarevich's conjecture}

\author{S\'andor J Kov\'acs \and Max Lieblich} 

\address{University of Washington, Department of Mathematics, Box 354350, Seattle, WA
  98195-4350, USA} 
\email{kovacs@math.washington.edu}
\email{lieblich@math.washington.edu}

\date{\today}

\thanks{S\'andor Kov\'acs was supported in part by NSF Grants DMS-0092165,
  DMS-0554697, and DMS-0856185, a Sloan Research Fellowship, and the Craig McKibben
  and Sarah Merner Endowed Professorship in Mathematics.  \newline %
  \indent Max Lieblich was supported in part by an NSF postdoctoral fellowship and by
  NSF Grant DMS-0758391.}

\subjclass[2000]{14J10, 14D23}

\begin{document}

\begin{abstract} We show that the number of deformation types of
  canonically polarized manifolds over an arbitrary variety with
  proper singular locus is finite, and that this number is uniformly
  bounded in any finite type family of base varieties. As a corollary
  we show that a direct generalization of the geometric version of
  Shafarevich's original conjecture holds for infinitesimally rigid
  families of canonically polarized varieties.
\end{abstract}

\maketitle
\tableofcontents
\markboth{S\'ANDOR J KOV\'ACS AND MAX LIEBLICH}{BOUNDEDNESS OF
  FAMILIES OF CANONICALLY POLARIZED MANIFOLDS}



\section{Introduction}
Fix an algebraically closed field $k$ of characteristic $0$. Let $B$ be a smooth
projective curve of genus $g$ over $k$ and $\Delta\subset B$ a finite subset.  A flat
morphism with connected fibers will be called a \emph{family}. For two families over
the same base $B$ a \emph{morphism of families} is simply a morphism of $B$-schemes.
A family $f: X\to B$ is called \emph{isotrivial} if $X_a\iso X_b$ for general points
$a,b\in B$, and $f: X\to B$ is \emph{admissible} (with respect to $(B, \Delta)$) if
it is not isotrivial and the map $f: X\setminus f^{-1}(\Delta)\to B\setminus \Delta$
is smooth.

At the 1962 International Congress of Mathematicians in Stockholm,
Shafarevich conjectured 
the following.

\begin{proclaimspecialr}{Shafarevich's Conjecture}
\label{shaff}
Let $(B, \Delta)$ be fixed and $q\geq 2$ an integer. Then
\begin{enumerate}
\item There exist only finitely many isomorphism classes of admissible
  families of curves of genus~$q$. \label{shafegy}
\item If\/ $2g-2+\#\Delta\leq 0$, then there exist no such
  families. \label{hyper} 
\end{enumerate}
\end{proclaimspecialr}

Shafarevich showed a special case of (\ref{shaff}.2): There exist no
smooth families of curves of genus $q$ over $\bP^1$. Conjecture
(\ref{shaff}) was proven by Parshin \cite{Parshin68} for
$\Delta=\emptyset$ and Arakelov \cite{Arakelov71} in general.

This conjecture has a natural analogue for curves over number fields.
For a brief discussion see \cite[\S 2]{Kovacs03c} and for more details
\cite{MR861969} and \cite{MR886677}.  Shafarevich's conjecture implies
Mordell's conjecture in both the function field and the number field
case by an argument known as Parshin's covering trick. Because of this
the proof of Shafarevich's conjecture in the number field case
constitutes the lion's share \cite{Aesop} of Faltings' celebrated
proof of Mordell's conjecture \cite{Faltings83,Faltings84}.

With regard to Shafarevich's conjecture, Parshin made the following observation. In
order to prove that there are only finitely many admissible families, one may proceed
as follows.  Instead of aiming for the general statement, first prove that there are
only finitely many deformation types\footnote{see~Definition~\ref{D:defo types}}. The
next step then is to prove that admissible families are rigid, that is, they do not
admit non-trivial deformations over a fixed base.
Now since every deformation type contains only one family, and since
there are only finitely many deformation types, the original statement
follows.

Based on this idea, the following reformulation of Shafarevich's
conjecture was used by Parshin and Arakelov to confirm the conjecture:

\begin{proclaimspecialr}{Shafarevich's Conjecture (version two)}
  Let $(B, \Delta)$ be fixed and $q\geq 2$ an integer. Then the
  following statements hold.
\begin{enumerate}
\item[$(\mathbf B)$] {\sc (Boundedness)} There exist only finitely
  many deformation types of admissible families of curves of genus $q$
  with respect to $B\setminus \Delta$.
\item[{$(\mathbf R)$}] {\sc (Rigidity)} There exist no non-trivial
  deformations of admissible families of curves of genus $q$ with
  respect to $B\setminus \Delta$.
\item[{$(\mathbf H)$}] {\sc (Hyperbolicity)} If $2g-2+\#\Delta\leq 0$,
  then no admissible families of curves of genus $q$ exist with
  respect to $B\setminus \Delta$.
\end{enumerate}
\end{proclaimspecialr}

\begin{rem}
  As we discussed above, $(\mathbf B)$ and $(\mathbf R)$ together
  imply (\ref{shaff}.1) and $(\mathbf H)$ is clearly equivalent to
  (\ref{shaff}.2).
\end{rem}

It is a natural and important question whether similar statements hold
for families of higher dimensional varieties. It is easy to see that
$(\mathbf R)$ fails \cite{Viehweg01}, \cite[10.4]{Kovacs03c} and hence
\eqref{shaff} fails. This gives additional importance to the
Parshin-Arakelov reformulation as it separates the clearly false part
from the rest. In fact, the past decade has seen a flood of results
concerning both $(\mathbf B)$ and $(\mathbf H)$.  For a detailed
historical overview and references to related results we refer the
reader to the survey articles \cite{Viehweg01}, \cite{Kovacs03c},
\cite{MVZ05}, \cite{Kovacs06a}, and \cite{Kovacs06b}.

In this article we are interested in $(\mathbf B)$. %
If there existed an algebraic stack $\mathfrak D$ parametrizing
families of canonically polarized varieties over the base $B\setminus
\Delta$, and if furthermore $\mathfrak D$ is of finite type, then
boundedness, $(\mathbf B)$, would follow. Before further discussing
the potential existence and properties of $\mathfrak D$, it behooves
us to mention the following notion closely related to $(\mathbf{B})$:

\begin{itemize}
\item[{$(\mathbf{WB})$}] \it {\sc (Weak Boundedness)} %
  We say that \emph{weak boundedness} holds if for an admissible family of projective
  varieties, $f\col X\to B$, the degree of $f_*\omega_{X/B}^m$ is bounded above in
  terms of $g(B), \#\Delta$, $m$, and $h_{X_{\rm gen}}$, where $X_{\rm gen}$ denotes
  the general fiber of $f$ and $h_{X_{\rm gen}}$ the Hilbert polynomial of
  $\omega_{X_{\rm gen}}^m$. In particular, the bound is independent of~$f$.
\end{itemize}

This was proven by Bedulev and Viehweg in 2000
\cite{Bedulev-Viehweg00}.  From this they derived the consequence that
as soon as a reasonable moduli theory exists for canonically polarized
varieties and if a $\mathfrak D$ as above exists, then it is indeed of
finite type.  Unfortunately, such $\mathfrak D$ almost never exist
(especially over open bases); moreover, when the base variety has
dimension higher than $1$, the question of how to rectify this
situation (by adding elements to the family over the discriminant
locus $\Delta$) is quite subtle.  The bulk of this paper is devoted to
pointing out that a proxy for $\mathfrak D$ can be constructed by
standard stack-theoretic methods, thus allowing us to show that
$(\mathbf{WB})$ implies $(\mathbf{B})$ while skirting the difficult
issues surrounding compactifications of the stack of canonically
polarized manifolds.

Before stating our main result we need the following definition.

\begin{defn} \label{D:defo types} %
  Let $U$ be a variety over a field $k$ and $\frC$ a class of schemes.  A morphism
  $X\to U$ is a $\frC$-morphism if for all geometric points $u\to U$, $X_u$ belongs
  to $\frC$. Two proper, flat $\frC$-morphisms $X_{1}\to U$, $X_{2}\to U$ are
  \emph{deformation equivalent\/} if there is a connected scheme $T$ with two points
  $t_{1},t_{2}\in T(k)$ and a proper, flat $\frC$-morphism $\mc X\to U\times T$ such
  that $\mc X|_{U\times t_i}\iso_{U}X_{i}$.  An equivalence class (with respect to
  deformation equivalence) of proper, flat $\frC$-morphisms $X\to U$ will be called a
  \emph{deformation type\/}.
\end{defn}

\begin{rem}
  In the sequel, the class $\frC$ will be chosen to be the class of
  canonically polarized varieties over a field $k$.
\end{rem}

The following theorem proves {$(\mathbf{B})$} in arbitrary dimension.

\begin{thm}\label{T:benjamin}
  Let $U$ be a variety over $k$ that is smooth at infinity (see
  \eqref{D:smooth at infinity}).  The set ${\rm Defo}_h(U)$ of
  deformation types of families $\mc X\to U$ of canonically polarized
  manifolds with Hilbert polynomial $h$ is finite.  Furthermore, if
  $T$ is a quasi-compact quasi-separated $\Q$-scheme and $\mc U\to T$
  is smooth at infinity, then there is an integer $N$ such that for
  every geometric point $t\to T$, we have 
  $\vert{\rm Defo}_h(\mc U_t)\vert\leq N$.
\end{thm}
This solves one of the open problems on the list compiled at the
American Institute of Mathematics workshop ``Compact moduli spaces and
birational geometry'' in December, 2004 \cite[Problem 2.4]{AIM04}. See
also \cite{Viehweg06a}.

In fact, we prove a more general result.  For the relevant terminology
see \S\ref{SS:main}.

\begin{thm}\label{T:main theorem}
  Let $\ms M^\circ$ be a weakly bounded\footnote{see
    Definition~\ref{D:weakly-bdd}} compactifiable Deligne-Mumford stack over a
  quasi-compact quasi-separated $\Q$-scheme $T$.  Given a morphism $U\to T$ that is
  smooth at infinity, there exists an integer $N$ such that for every geometric point
  $t\to T$, the number of deformation types of morphisms $U_t\to\ms M^\circ_t$ is
  finite and bounded above by $N$.
\end{thm}


\begin{demo}{\bf Definitions and Notation}\label{demo:defs}
  
  For morphisms $f:X\to B$ and $\vartheta: T\to B$, the symbol $X_T$ will denote
  $X\times_B T$ and $f_T:X_T\to T$ the induced morphism. Of course, by symmetry we
  also have the notation $\vartheta_X:T_X\simeq X_T\to X$. In particular, for $b\in
  B$ we write $X_b = f^{-1}(b)$. In addition, if $T=\Spec F$, then $X_T$ will also be
  denoted by~$X_F$.
  Finally, if $\sF$ is an $\sO_X$-module, then $\sF_T$ will denote the
  $\sO_{X_T}$-module $\vartheta_X^*\sF$.

  Given a proper scheme $X$ over a field $k$, we write $\Pic^{\tau}_X$
  for the locus of numerically trivial invertible sheaves in $\Pic_X$.
  This is generally larger than $\Pic^0_X$, the connected component containing the
  trivial sheaf.  Given a field extension $L/k$ and an invertible sheaf $\ms N$ on
  $X_L$, we will write $[\ms N]$ for the element of $\Pic_X(L)$ associated to $\ms
  N$.

  For the theory of stacks, we will use the definitions and
  conventions of \cite{Lau-MB00}.  In particular, all algebraic
  (Deligne-Mumford or Artin) stacks are assumed to be quasi-separated.
  Most of the time, the stacks we use will in fact be separated; this
  is always indicated in the text as a hypothesis when it is used.

  A quasi-compact separated Deligne-Mumford stack $\ms M$ is
  \emph{polarized\/} if there exists an invertible sheaf $\ms L$ on
  $\ms M$ such that the non-vanishing loci of all sections of all
  tensor powers of $\ms L$ generate the topology on the underlying
  topological space of $\ms M$.  Equivalently, some tensor power of
  $\ms L$ is the pullback from the coarse moduli space $M$ of an ample
  invertible sheaf $L$.  A Deligne-Mumford stack is \emph{tame\/} if
  the order of the stabilizer group of any geometric point $x$ is
  invertible in $\kappa(x)$.  We will only explicitly encounter tame stacks in
  the generalities of Section \ref{S:stable}; otherwise, we will be
  working in characteristic $0$, where tameness is automatic and will
  go unmentioned.
  %
\end{demo}

\begin{demo*}{\bf Acknowledgments}
  Our work started after an inspiring conversation with Eckart Viehweg at the AIM
  workshop mentioned above.  We have also benefited from discussions with Brian
  Conrad, Brendan Hassett, J\'anos Koll\'ar, Martin Olsson, and Paul Smith. We would
  like to thank them all for their helpful comments.

  We are grateful to AIM for catalyzing our research collaboration and would like to
  thank the referees for several lists of useful comments, corrections, and
  suggestions.
\end{demo*}

\section{Coarse boundedness}

\subsection{Bounding maps to a projective scheme}

As in the introduction, $k$ will be an algebraically closed field of
characteristic $0$. In what follows, \emph{variety} will mean a
\emph{$k$-variety}.

\begin{defn} \label{D:smooth at infinity} For a morphism $\mc U\to T$ of algebraic
  spaces let $\Sing(\mc U/T)$ denote the smallest closed subset of $\mc U$ such that
  the induced morphism $\mc U\setminus \Sing(\mc U/T)\to T$ is smooth. The morphism
  $\mc U\to T$ is called \emph{smooth at infinity\/} if it is of finite presentation,
  $\Sing(\mc U/T)$ is proper over $T$, and $\mc U\setminus \Sing(\mc U/T)$ is
  schematically dense in every geometric fiber.  A variety will be called smooth at
  infinity if its structure morphism is smooth at infinity.
\end{defn}


\begin{newnumr}\label{Para:initial fun}
  Let $M$ be a proper $k$-scheme with a fixed invertible sheaf $\ms N$
  and let $U$ be an algebraic variety that is smooth at infinity.  By
  Nagata's theorem, $U$ embeds into a proper variety $B$.  Blowing up
  and using the assumption that $U$ is smooth at infinity, and hence
  $\Sing U$ is proper, we may assume that $B\setminus U$ is a divisor
  $\Delta$ (with simple normal crossings, if desired) and that $B$ is
  smooth in a neighborhood of $\Delta$. Because $M$ is proper, it
  follows that given a morphism $\xi:U\to M$, there is an open subset
  $\iota:U'\into B$ containing $U$ and every codimension $1$ point of
  $B$ and an extension of $\xi$ to a morphism $\xi':U'\to M$. Taking
  the reflexive hull of $\iota_*\ms N_{U'}$ yields an invertible sheaf
  $\ms N_{\xi}$ on $B$ by \cite[II.1.1.15, p.154]{OSS}.
  
  On the other hand, suppose $C^\circ$ is a smooth curve over $k$ with
  smooth compactification $C$.  Given a morphism $C^\circ\to U$ and a
  morphism $\xi:U\to M$ as above, one obtains an extension $\xi_C:C\to
  M$ of the restriction of $\xi$ to $C^\circ$.  It is of course not
  necessary for $\deg(\xi_C^*\ms N)$ to equal $\deg\ms N_{\xi}|_C$,
  but this will clearly occur when $C$ is contained in $U'$ (in the
  above notation).
\end{newnumr}

\begin{defn}
  A \emph{$(g,d)$-curve\/} is a smooth curve $C^\circ$ whose smooth
  compactification $C$ has genus $g$ and such that $C \setminus
  C^\circ$ consists of $d$ closed points.
\end{defn}

\begin{defn}\label{D:coarse weak boundedness} Given $U$ and $M$ as
  above, a morphism $\xi:U\to M$ is \emph{weakly bounded with respect
    to $\ms N$\/} if there exists a function $\wb_{\ms N}:\Z_{\geq
    0}^2\to\Z$ such that for every pair $(g,d)$ of non-negative
  integers, for every $(g,d)$-curve $C^\circ\subseteq C$, and for
  every morphism $C^\circ\to U$, one has that $\deg \xi_C^*\ms N\leq
  \wb_{\ms N}(g,d)$.  The function $\wb_{\ms N}$ will be called a
  \emph{weak bound\/} (with respect to $\ms N$), and we will say that
  \emph{$\xi$ is weakly bounded by $\wb_{\ms N}$\/}.
\end{defn}

\begin{notn} Given a field extension $L/k$, the set of morphisms $U_L\to
  M_L$ which are weakly bounded by $\wb_{\ms N}$ will be denoted
  $\setW(U,M,\wb_{\ms N})(L)$. Notice that as $\wb_{\sN}$ depends on
  $\sN$, so does $\setW(U,M,\wb_{\ms N})(L)$.
\end{notn}

\begin{prop}\label{P:small family} Let $\wb$ be a weak bound. Then
  there exists a variety $\WW^{\wb}$ and a morphism
  $\Xi:\WW^{\wb}\times U\to M$ such that for every field extension
  $L/k$ and for every morphism $\xi:U_L\to M_L$ that is weakly bounded
  by $\wb$ there exists an $L$-valued point $p:\Spec L\to \WW^{\wb}$
  such that $\xi=\Xi\resto{\{p\}\times U}$.
\end{prop}


\begin{rem}
  Notice that this does not necessarily mean that every point of $\WW^{\wb}$
  corresponds to a weakly bounded morphism $U\to M$.  This phenomenon is common in
  the theory of moduli; one often produces a bounded family containing the points of
  interest, but possibly also containing numerous other points. 
  In fact, this is one of the main difficulties in the present situation.  It is much easier
  to find a bounding family than one that actually parametrizes the class we are
  interested in.
  %
\end{rem}

The proof consists of several steps.  First, we compactify $U\subseteq
B$ as the complement of a divisor $\Delta$ in a proper variety as in
\eqref{Para:initial fun}.  Then we bound the set of invertible sheaves
$\ms N_{\xi}$. The choice of $n+1$ sections of such an $\ms N_{\xi}$
that simultaneously vanish only in $\Delta$ can then be parametrized
by a finite type space $T$.

\begin{demo}{Assumption}\label{ass-M-is_proj}
  We will assume that $M$ is projective, fix an embedding $M\into \P^n$, and let $\ms
  N=\sO_M(1)=\ms O_{\P^n}(1)\resto M$. For simplicity we replace the phrase ``weakly
  bounded with respect to $\ms N$'' by ``weakly bounded''.
\end{demo}

Let us first treat the case $M=\bP^n$.

\begin{lem}\label{lem:constr-W-in-Pic} 
  Given a compactification $U\subseteq B$ as above, there exists a reduced
  subscheme of finite type $\schW(U,\P^n,\wb)\subset\Pic_{B}$ such that for all field
  extensions $L/k$ and for all $\xi\in\setW(U,\P^n,\wb)(L)$, we have $[\ms
  N_{\xi}]\in\schW(U,\P^n,\wb)(L)$.
\end{lem}

\begin{proof} We first claim that it suffices to prove the result when $B$ is smooth
  and projective.  Indeed, choose a projective resolution of singularities (or a
  projective alteration \cite{MR1423020}) $\pi:\widetilde B\to B$ with $\widetilde B$
  smooth, and let $\widetilde U$ be the preimage of $U$.  Now we can consider weakly
  bounded morphisms $\widetilde U\to M$ with the same weak bound $\wb$.  Among these
  will be the compositions of $\pi$ with morphisms $U\to M$ weakly bounded by $\wb$.
  In other words, composition with $\pi$ induces a natural map $\pi^*:
  \setW(U,\P^n,\wb)\to \setW(\widetilde U,\P^n,\wb)$. Observe that the pullback
  morphism $\pi^{\ast}:\Pic_B\to\Pic_{\widetilde B}$ is of finite type by
  \cite[XII.1.1]{SGA6}, and hence if the required $\schW(\widetilde
  U,\P^n,\wb)\subset\Pic_{\widetilde B}$ exists, then
  $\schW(U,\P^n,\wb):=(\pi^\ast)^{-1}\schW(\widetilde U,\P^n,\wb)\subset \Pic_B$
  satisfies the desired conditions. Therefore we may assume from now on that $B$ is
  smooth and projective.

  Suppose $\dim B\geq 3$ and let $Y\subset B$ be a general ample divisor.  By
  \cite[XIII.3.8]{SGA6}, the restriction morphism $\Pic_B\to\Pic_Y$ of Picard schemes
  is of finite type.  Since the restriction of a morphism $U\to\P^n$ weakly bounded
  by $\wb$ to $U\cap Y$ is also weakly bounded by $\wb$ and we have $\ms N_\xi\resto
  Y\iso \ms N_{\xi_Y}$, we see that it suffices to prove the statement for $U\cap
  B\subset B$.  Thus, we may assume $\dim B\leq 2$.

  If $\dim B=1$, then the inclusion $U\subset B$ is a $(g,d)$-curve with $g$ the
  genus of $B$ and $d$ the number of points in $B\setminus U$.  In addition, any
  morphism $\xi:U\to \P^n$ extends to a morphism $B\to\P^n$.  By the weak boundedness
  assumption, $0\leq\deg\xi_C^{\ast}\ms N\leq\wb(g,d)$, so that $\xi_C^{\ast}\ms N$
  is contained in the preimage $\Pic_C^{[0,\wb(g,d)]}$ of the interval
  $[0,\wb(g,d)]\subset\Z$ under the degree map $\Pic_C\to\Z$.  Since the fiber of
  $\Pic_C\to\Z$ over any finite subset is of finite type, we see that setting
  $\schW(U,\P^n,\wb)=\Pic_C^{[0,\wb(g,d)]}$ yields the result.

  Hence we may assume for the rest of the proof that $B$ is a surface.  Let $A$ be a
  very ample divisor on $B$. We will prove that for each $\xi$ we have
  $0\leq\deg_A\ms N_{\xi}\leq N$ and $c_1(\ms N_\xi)^2\geq 0$.  These conditions
  define an open subscheme $\WW(U,\P^n,\wb)$ of $\Pic_B$. Moreover, by \cite[Theorem
  XIII.3.13(iii)]{SGA6}, there exists a quasi-compact scheme $T$ and a family of
  invertible sheaves $\sL$ on $B\times T$ such that every sheaf $\ms N_{\xi}$ appears
  as a fiber over a point of $t$.  We conclude that $\WW(U,\P^n,\wb)$ is
  quasi-compact and therefore of finite type, as desired.

  So it remains to verify that the above numerical conditions are satisfied.  We may
  assume that the very ample divisor $A$ is smooth. The definition of weak
  boundedness then yields a bound $\deg\xi_A^{\ast}\ms O_M(1)\leq N$ which depends
  only on the genus of $A$ and on $A\cdot\Delta$.  Moreover, since $\codim(B\setminus
  U', B)\geq 2$, we can choose an $A$ such that $A\subset U'$. In this case $\ms
  N_{\xi}|_A\iso \xi_A^{\ast}\ms O_M(1)$ and hence we conclude that $0\leq\deg_A\ms
  N_{\xi}\leq N$.
  
  Next consider $H_1, H_2$, the zero loci of two general sections of $\ms O_M(m)$ for
  some $m\gg 0$. Assume that $\xi$ is non-constant and let $\widehat
  H_i=\widebar{\xi'^*H_i}\subset B$ be the closure of the pullback of $H_i$ to $U'$
  via $\xi'$ for $i=1,2$.  Clearly, $\widehat H_1\cdot \widehat H_2\geq 0$. Notice that by
  definition $\ms N_{\xi}\iso \ms O_B(\widehat H_i)$ for $i=1,2$ and hence $c_1(\ms
  N_\xi)^2\geq 0$, as desired.
\end{proof}

\begin{lem}\label{L:sections} 
  Given a finite type reduced subscheme $Y\subset\Pic_{B}$, there is a
  finite stratification $Y_i$ of $Y$ by locally closed subschemes such
  that the functor of tuples $(\ms L,\sigma^0,\dots,\sigma^{n})$ with
  $[\ms L]\in\coprod Y_i$ and $\sigma^0,\dots,\sigma^{n}$ 
  global sections of $\ms L$ at least one of which is non-zero is represented by a
  reduced scheme $W$ separated and of finite type over $\coprod Y_i$.
\end{lem}
\begin{proof} 
Let $\ms Pic_B$ be the Artin stack of invertible sheaves on $B$; this
is a $\G_m$-gerbe over the Picard scheme $\Pic_B$.  Write $\ms P\to Y$
for the fiber product $\ms Pic_B\times_{\Pic_B}Y$.  Write $\ms
L_{\text{univ}}$ for the universal invertible sheaf on $\ms P\times
B$.  The function $\ms L\mapsto\operatorname{h}^0(\ms L)$ is
upper semicontinuous and thus defines a stratification of $\ms P$ by
reduced locally closed substacks $\ms P_i\subset\ms P$.  Since $\ms
P\to Y$ is a $\G_m$-gerbe, each $\ms P_i$ is a $\G_m$-gerbe over a
locally closed subscheme $Y_i\subset Y$, and the $Y_i$ form a
stratification of $Y$.  

Write $p_i:\ms P_i\times B\to\ms P_i$ for the first projection. By
cohomology and base change, the sheaf $(p_i)_\ast\ms
L_{\text{univ}}|_{\ms P_i\times B}$ is a locally free $\ms
P_i$-twisted sheaf (see \cite[Section 3.1.1]{MR2388554} for the
definition and basic properties of twisted sheaves).  A choice of sections
$\sigma^0,\ldots,\sigma^n$ such that at least one $\sigma^i$ is not
the zero section is a point of the stack 
$$
W:=\V((((p_i)_\ast\ms L_{\text{univ}}|_{\ms P_i\times
  B})^{\vee})^{\oplus(n+1)})\setminus \mathbf 0,
$$ where $\mathbf 0$ denotes the zero
section of the vector bundle.

Since the inertia stack of $\ms P$ acts on $\ms L_{\text{univ}}$ by scalar
multiplication, the induced action on $W$ is faithful, from which it follows that $W$
is an algebraic space.  To prove that it is a separated scheme of finite type (and to
give a more concrete description of the space), we can work \'etale locally on the
sheafification $P_i$ of $\ms P_i$ and thus assume that (1) $\ms P_i$ is isomorphic to
$P_i\times\B\G_m$ and (2) the pullback of $(p_i)_\ast\ms L_{\text{univ}}$ via the
canonical map $P_i\to P_i\times\B\G_m$ is trivial, say of rank $r$.  The natural left
action of $\G_m$ on the fibers of the locally free sheaf $(((p_i)_\ast\ms
L_{\text{univ}}|_{\ms P_i\times B})^{\vee})^{\oplus n}$ is via scalar multiplication.
Thus, $W$ is isomorphic to the stack-theoretic quotient of the scalar multiplication
action of $\G_m$ on $\A_{P_i}^{nr}\setminus\mathbf 0$.  This is just
$\P^{nr-1}_{P_i}$, which is certainly separated and of finite type.  (Continuing
along these lines shows that $W$ is in fact isomorphic to a Brauer-Severi scheme over
$P_i$ with the same Brauer class as $[\ms P_i]$.  While it may seem baffling that a
projective $P_i$-scheme can be an open substack of a geometric vector bundle over
$\ms P_i$, it arises from the fact that $\ms P_i$ --- and therefore any vector bundle
over $\ms P_i$ --- is highly non-separated.)

\end{proof}

\begin{lem}\label{L:inclusion} 
  Let $S$ be a reduced Noetherian algebraic space and $\ms X\to S$ and
  $\ms Y\to S$ two Artin stacks of finite presentation. Let $\ms
  Z\subset\ms X$ and $\ms T\subset\ms Y$ be locally closed
  substacks. Given an $S$-morphism $\phi:\ms X\to\ms Y$, there is a
  monomorphism of finite type $S'\to S$ whose image contains a
  geometric point $s\to S$ if and only if $\phi_{s}$ maps $(\ms
  Z_{s})_{\text{red}}$ into $(\ms T_{s})_{\text{red}}$.
\end{lem}
\begin{proof}
  Pulling back $\ms T$ to $\ms X$ and replacing the inclusion $\ms
  T\subset\ms Y$ by $\ms T_{\ms X}\subset\ms X$, we may assume that
  $\ms X\to\ms Y$ is the identity morphism.  Then the set $\ms
  Z\setminus\ms T$ is constructible in $\ms X$\!, so the reduced
  structure on the complement of its image in $S$ is constructible.
  Any constructible set admits a natural locally finite stratification
  by reduced algebraic spaces, yielding the desired morphism $S'\to
  S$.
\end{proof}


\begin{lem} \label{lem:constructing-W-Delta}
  Let $W$ be the scheme constructed in \eqref{L:sections}.  Then there
  is a finite type morphism $W_{\Delta}\to W$ such that for any $w\in
  W$ the reduced common zero locus of $\sigma^0_w,\cdots,\sigma^{n}_w$
  is contained in $\Delta$ if and only if $w$ is in the image of
  $W_{\Delta}$.  In fact, $W_{\Delta}$ is the union of pieces in a
  stratification of $W$.
\end{lem}
\begin{proof} The sections $\sigma^i$ define divisors $Z_i\subset
  B\times W$. Apply \eqref{L:inclusion} with $\sX=\sY=B\times W$, $S=W$,
  $\sZ=Z_0\cap\cdots\cap Z_{n}$ and $\sT=\Delta\times W$.
\end{proof}

\begin{proof}[Proof of \eqref{P:small family}] %
  Let $\WW^{\wb}=W_{\Delta}$ be the result of applying \eqref{L:sections} and
  \eqref{lem:constructing-W-Delta} to the scheme $Y=\WW(U,\bP^n,\wb)$ constructed in
  \eqref{lem:constr-W-in-Pic}.  The sections $\sigma^0,\ldots,\sigma^n$ define the
  required morphism $\Xi: \WW^{\wb}\times U\to\P^n$, proving the statement for
  $M=\P^n$.

  For a general projective $M\inj\P^n$, if we let $\WW^{\wb}\times
  U\to\WW^{\wb}\times \P^n$ be the morphism ensured by the previous case, we can take
  $S=\WW^{\wb}$, $\ms Z=\ms X=\WW^{\wb}\times U$, $\ms Y=\WW^{\wb}\times \P^n$, and
  $\ms T=\WW^{\wb}\times M$ in \eqref{L:inclusion}, yielding a finite type
  monomorphism $W'\to\WW^{\wb}$ and the required morphism $W'\times U\to M$.
\end{proof}

\begin{prop} \label{P:2.7} Given a polarized variety $(M, \ms O_M(1))$, an open
  subscheme $M^\circ$ and a weak bound $\wb$, there is a $k$-variety
  $\WW_{M^\circ}^{\wb}$ and a morphism $\WW_{M^\circ}^{\wb}\times U\to M^\circ$ such
  that for every field extension $L/k$, every morphism $U_L\to M_L^\circ$ whose
  composition with the inclusion $M_L^\circ\into M_L$ is weakly bounded with
  respect to the polarization of $M_L$ by $\wb$ appears in a fiber over
  $\WW_{M^\circ}^{\wb}(L)$.
\end{prop}
\begin{proof}
  This follows from \eqref{P:small family} and \eqref{L:inclusion}.
\end{proof}

We briefly indicate how to extend the results above to the case of a
family over a reduced base.
\begin{prop} \label{P:coarse-compactification} Let $T$ be a
  quasi-compact quasi-separated reduced scheme and $\mc U\to T$ a
  separated morphism which is smooth at infinity.  Given a projective $T$-scheme
  of finite presentation $(\mc M, \ms O_{\mc M}(1))$, an open subscheme $\mc
  M^\circ\subseteq \mc M$ of finite presentation over $T$, and a weak
  bound $\wb$, there exists a $T$-scheme of finite presentation
  $\WWW_{\mc M^\circ}^{\wb}$ and a morphism $\Xi:\WWW_{\mc
    M^\circ}^{\wb}\times \mc U\to \mc M^\circ$ such that for every
  geometric point $t\to T$ and for every morphism $\xi:\mc U_t\to \mc
  M_t^\circ\subseteq \mc M_t$ that is weakly bounded by $\wb$ there
  exists a point $p\to \WWW_{\mc M_t^\circ}^{\wb}$ such that
  $\xi=\Xi\resto{\{p\}\times \mc U_t}$.
\end{prop}

\begin{proof} Since $T$ is quasi-compact and quasi-separated, absolute
  Noetherian approximation \cite[C.3 and Theorem C.9]{MR1106918} lets
  us assume that $T$ is of finite type over $k$.  We claim that there
  is a finite type morphism $T'\to T$ and a fiberwise dense open
  immersion $U_{T'}\inj\mc B$ with $\mc B\to T'$ a proper scheme with
  geometrically integral fibers.  To see this, we can first replace
  $T$ with the disjoint union of its irreducible components and thus
  assume that $T$ is integral.  The geometric generic fiber
  $U_{\widebar\eta}$ has an integral compactification
  $U_{\widebar\eta}\inj\widebar B$ by Nagata's theorem.  Since
  $\widebar B$ is of finite presentation over $T$, there is a finite
  type integral $T$-scheme $T_1\to T$ with a lift $\widebar\eta\to
  T_1$ to a geometric generic point over the given geometric generic
  point of $T$, a proper $T_1$-scheme $\mc B\to T_1$, and an open
  immersion $U_{T_1}\inj\mc B$ whose pullback to $\widebar\eta$ is
  $U_{\widebar\eta}\inj\widebar B$.  

  By generic flatness, there is a dense open subscheme $T_2\subset T_1$ over which
  $\mc B$ and $\mc B_{T_2}\setminus U_{T_2}$ are flat.  Applying \cite[Theorem
  12.2.4]{EGAIV4}, there is a further open subscheme $T_3\subset T_2$ over which the
  geometric fibers of $\mc B$ are integral and the geometric fibers of $\mc
  B_{T_3}\setminus\mc U_{T_3}$ have dimension strictly smaller than the fiber
  dimension of $\mc B_{T_3}\to T_3$.  It follows that $U_{T_3}\subset\mc B_{T_3}$ is
  a fiberwise dense open immersion.  Since $T_3\to T$ is dominant and of finite type,
  Chevalley's theorem shows that its image contains a dense open.  By Noetherian
  induction, there are morphisms $T'\to T$ and $U_{T'}\inj\mc B$ as claimed in the
  first paragraph.  Resolution of singularities and similar stratification and
  Noetherian induction argument gives a finite type morphism $T''$ over which there
  is a smooth projective morphism $\widetilde{\mc B}\to T''$ whose geometric fibers
  are connected and which admits a fiberwise birational morphism $\widetilde{\mc
    B}\to\mc B_{T''}$.  We can replace $T$ by $T''$ and assume that we have such a
  compactification and resolution.

  Since the geometric fibers of $\mc B\to T$ are integral, we see that
  $\mc B$ is cohomologically flat in degree $0$ and the Picard functor
  is separated.  Thus, the Picard stack $\ms Pic_{\mc B/T}$ is a
  $\G_m$-gerbe over a separated algebraic space $\Pic_{\mc B/T}$
  locally of finite type over $T$.  Given a $T$-flat relatively ample smooth
  divisor $\mc B'\subset\mc B$, if the fiber dimension of $\mc B'$ is at
  least $2$ then by \cite[Theorem XIII.3.8]{SGA6} the restriction
  morphism $\ms Pic_{\mc B/T}\to\ms Pic_{\mc B'/T}$ is of finite
  type.  Since we can always replace $T$ by an open covering, we can
  always assume that such a $T$-smooth divisor exists.

  We claim that there is an open substack $\ms Y\subset\ms Pic_{\mc
    B/T}$ of finite type over $T$ such that for each geometric point
  $t\to T$, the invertible sheaves $\ms N_{\xi}$ on $B_t$ arising from
  weakly bounded morphisms $U_t\to\mc M^\circ$ via the procedure of
  \eqref{Para:initial fun} lie in $\ms Y$.  Arguing precisely as in
  the proof of \eqref{lem:constr-W-in-Pic}, it suffices to prove this
  for the inclusion $\widetilde U\inj\widetilde{\mc B}$, so that we
  may assume $\mc B$ is smooth over $T$.

  First suppose $\mc B$ is a relative smooth curve.  The degree map
  gives a finite type morphism $\Pic_{\mc B/T}\to\Z_T$.  Just as in
  the proof of \eqref{lem:constr-W-in-Pic}, the weak boundedness shows
  that the $\ms N_{\xi}$ are contained in the preimage of a finite
  interval in $\Z$, yielding the claim.

  If the fibers of $\mc B/T$ have dimension at least $2$, then arguing
  as in the proof of \eqref{lem:constr-W-in-Pic} and using the
  existence of smooth relatively ample divisors $\mc B'\subset\mc B$
  (after possibly replacing $T$ by an open covering), we can reduce to
  the case in which $\mc B\to T$ is a smooth projective relative
  surface.  Now, again as in the proof of \eqref{lem:constr-W-in-Pic}
  we have that the sheaves $\ms N_{\xi}$ on the fibers over $T$
  satisfy $0\leq\deg_A\ms N_{\xi}\leq N$ and $c_1(\ms N_{\xi})^2\geq
  0$ (where $A\subset\mc B$ is a relatively ample smooth divisor).
  Invoking \cite[Theorem XIII.3.13(iii)]{SGA6} again, we see that the
  collection of invertible sheaves on the fibers satisfying those
  boundedness conditions forms a finite type open substack $\ms
  Y\subset\ms Pic_{\ms B/T}$, as desired.

  Now, to bound the map we argue as in (\ref{L:sections}).  As
  written, the argument is completely general and applies in the
  present situation.  It yields the universal collection of sections.
  The proofs of (\ref{L:inclusion}), (\ref{lem:constructing-W-Delta}),
  and (\ref{P:2.7}) also carry over to yield the map $\Xi$.
\end{proof}

\begin{newnumr}
Next we compactify $\xi$ in a bounded family.  
\end{newnumr}

\begin{defn} Given a $T$-scheme $B\to T$, a \emph{relative \sncd\/}
  $D\subset B$ is a divisor of the form $D=D_1+\cdots+D_r$ such that
  $B$ is flat over $T$ in a neighborhood of $D$, each $D_i$ is flat
  over $T$, and in each geometric fiber $B_{t}$ the divisor
  $(D_1)_{t}+\cdots+(D_r)_{t}$ is a \sncd.
\end{defn}

\begin{prop}\label{L:compactifying family} Let $T$ be reduced
  quasi-separated and quasi-compact and $\mc U\to T$ a separated
  morphism that is smooth at infinity.  Given a proper $T$-scheme of
  finite presentation $\mc M$ and a $T$-morphism $\xi:\mc U\to \mc M$,
  there exists a finite type surjective morphism $T'\to T$, a proper scheme $\mc
  B\to T'$, an open immersion $\mc U_{T'}\inj \mc B$ over $T'$ whose
  complement $\mc B\setminus \mc U_{T'}$ is a relative \sncd, and a
  $T'$-morphism $\widebar{\xi}:\mc B\to \mc M_{T'}$ such that
  $\widebar{\xi}|_{\mc U_{T'}}=\xi_{T'}$.
\end{prop}

\begin{proof} 
  By absolute Noetherian approximation, we may assume that $T$ is Noetherian.  We may
  then replace $T$ by the disjoint union of its irreducible components and assume
  that $T$ is integral. Next we compactify the morphism $\mc U\to T$ to a proper
  scheme $\mc B'\to T$ (which is not necessarily flat!).  Resolving the singularities
  of the generic fiber of $\mc B'\setminus{\Sing(\mc U/T)}$ yields an immersion $\mc
  U\to \mc B''$ into a proper scheme over the function field of $T$ whose general
  fiber over $T$ is smooth outside $\mc U$.  After a birational modification of $\mc
  B''$, we may assume that $\xi$ extends to $\mc B''$ and that $\mc B''\setminus \mc
  U$ is a \sncd.  This extends over an open dense subscheme of $T$.  By Noetherian
  induction, we can thus stratify $T$ so that such compactifications exist over each
  stratum.  Given the compactifications, we proceed as in the proof of \eqref{P:small
    family}.
\end{proof}

\begin{remark}\label{hypothesis-whore}
  If $\mc U\to T$ is quasi-projective, then we can assume that $\mc
  B\to T'$ is projective, as the resolution of singularities of $\mc
  B'$ can be assumed to be projective.
\end{remark}

\subsection{Bounding maps to a quasi-projective scheme}

\begin{defn}\label{D:rel-ample-wrt}
  Given a proper $T$-scheme $\pi:\mc M\to T$ and an open subscheme
  $\mc M^\circ\subseteq \mc M$, an invertible sheaf $\ms L$ on $\mc M$
  is \emph{relatively ample with respect to $\mc M^\circ$} if there
  exists an integer $m>0$ such that
  \begin{enumerate}
  \item $\pi^*\pi_* \ms L^m \to \ms L^m$ is surjective over $\mc
    M^\circ$, and
  \item the natural map $\mc M^\circ\rightarrow \P_T \left( \pi_*\ms
      L^m\right)$ is a locally closed immersion.
  \end{enumerate}
  A \emph{relative polarization of $\mc M$ with respect to $\mc
    M^\circ$} is an invertible sheaf $\ms L$ that is relatively ample
  with respect to $\mc M^\circ$.
\end{defn}

Note that if (\ref{D:rel-ample-wrt}.2) holds for some $m>0$, then it
holds for any $m$ sufficiently large and divisible.

\begin{defn}
  Given a separated $T$-scheme of finite type $\mc M^\circ$, a
  \emph{relative compactification of $\mc M^\circ$} is a $T$-morphism
  $\iota: \mc M^\circ\to \mc M$ that embeds $\mc M^\circ$ as an open
  subscheme of the proper $T$-scheme $\mc M$. If there is no danger of
  confusion, we will abuse notation and refer to a relative
  compactification $\iota: \mc M^\circ\to \mc M$ simply as $\mc M$.  A
  \emph{morphism} between relative compactifications $\iota: \mc
  M^\circ\to \mc M$ and $\iota': \mc M^\circ\to \mc M'$ is a
  $T$-morphism $\phi: \mc M\to \mc M'$ such that
  $\phi\circ\iota=\iota'$.
\end{defn}

\begin{rem}
  These notions seem most natural if $\mc M^\circ_t$ is dense in $\mc
  M_t$ for all $t\in T$, but we do not need to make this assumption
  here.
\end{rem}

The next statement allows us to replace a polarization with respect to an open
subscheme with a polarization everywhere.

\begin{prop}\label{P:love-triangle}
  Let $T$ be a Noetherian scheme, $\mc M^\circ$ a separated $T$-scheme
  of finite type, $\iota:\mc M^\circ\to \mc M$ a relative
  compactification, and $\ms L$ a relative polarization of $\mc M$
  with respect to $\mc M^\circ$. Then there exists a diagram of
  relative compactifications of $\mc M^\circ$
  $$
  \xymatrix{%
    & \widetilde{\mc M} \ar[dl]_\sigma \ar[dr]^\tau & \\
    \mc M & \mc M^\circ \ar[l]^\iota \ar[u]^{\tilde\iota}
    \ar[r]_{\iota'}& \mc M' }
  $$
  and a $T$-ample invertible sheaf $\sA$ on $\mc M'$ such that
  \begin{enumerate}
  \item there exists an inclusion of invertible sheaves
    $\tau^*\sA\subseteq \sigma^*\sL$ which is an isomorphism on $\mc
    M^\circ$.
  \item In particular, given a weak bound $\wb$, a geometric point
    $t\to T$, and a morphism $\xi : U\to \mc M^\circ_t$, if
    $\iota\circ\xi$ is weakly bounded with respect to $\ms L_t$ by
    $\wb$, then $\iota'\circ\xi$ is weakly bounded with respect to
    $\ms A_t$ by $\wb$.
  \end{enumerate}
\end{prop}

\begin{proof}
  Let $m> 0$ be the integer given in \eqref{D:rel-ample-wrt} and $\ms
  E=\pi_*\ms L^m$.  Consider the natural map $\nu:\pi^*\ms
  E=\pi^*\pi_*\ms L^m\to \ms L^m$ which is surjective on $\mc
  M^\circ$. Let $\ms I=\nu(\pi^*\ms E)\otimes \ms L^{-m}\subseteq \ms
  O_{\mc M}$ and let $\sigma:\widetilde{\mc M}\to \mc M$ be the
  blowing up of the ideal sheaf $\ms I$.

  Since the support of $\factor{\ms O_{\mc M}}.{\ms I}.$ is disjoint
  from $\mc M^\circ$, $\tilde\iota=\sigma^{-1}\circ\iota:\mc
  M^\circ\to \widetilde{\mc M}$ is a relative compactification of $\mc
  M^\circ$ and $\ms N=\sigma^*\ms L^m\otimes \sigma^{-1}\ms I \cdot\ms
  O_{\widetilde{\mc M}}$ is relatively ample with respect to
  $\tilde\iota(\mc M^\circ)$. The surjective morphism
  $\sigma^*\pi^*\ms E\onto \ms N$
  induces a $T$-morphism $\tau: \widetilde{\mc M}\to \P_T(\ms E)$
  which is an embedding on $\mc M^\circ$. Letting $\mc M'$ be the
  scheme-theoretic image of $\tau$ and $\ms A$ the restriction of $\ms
  O_{\P_T(\ms E)}(1)$ to ${\mc M'}$ yields (\ref{P:love-triangle}.1).

  Given a curve $C$ and a morphism $\gamma: C\to \widetilde{\mc M}_t$
  such that $\gamma(C)\cap \mc M^\circ_t\neq \emptyset$, the natural
  map $\gamma^*\tau^*\ms A\to \gamma^*\sigma^*\ms L$ remains
  injective. Therefore (\ref{P:love-triangle}.1) implies
  (\ref{P:love-triangle}.2).
\end{proof}

\begin{cor}\label{C:cc-quasi-proj}
  Let $T$ be a Noetherian scheme, $\mc M^\circ$ a separated $T$-scheme
  of finite type, $\iota:\mc M^\circ\to \mc M$ a relative
  compactification, $\ms L$ a relative polarization of $\mc M$ with
  respect to $\mc M^\circ$, and $\wb$ a weak bound. Then there exists
  a $T$-scheme of finite type $\WWW_{\mc M^\circ}^{\wb}$ and a
  morphism $\WWW_{\mc M^\circ}^{\wb}\times \mc U\to \mc M^\circ$ such
  that for every geometric point $t\to T$, every morphism $\mc U_t\to
  \mc M_t^\circ\subseteq \mc M_t$ which is weakly bounded with respect
  to $\ms L_t$ by $\wb$ appears in a fiber over $\WWW_{\mc
    M_t^\circ}^{\wb}$.
\end{cor}

\begin{proof}
  This follows directly from \eqref{P:coarse-compactification} and
  \eqref{P:love-triangle}.
\end{proof}



\section{Weak stacky stable reduction}\label{S:stable}

\subsection{Groupoid-equivariant objects in a stack}
We start with a few basic results about equivariant objects and their
liftings.  While the main result of this section can be stated in
purely stack-theoretic language (as we indicate in the alternative
proof of \eqref{C:fluffy lift}), the formalism we briefly sketch here
is useful for clarifying the proof of \eqref{P:extension}.  

Let $(R,Z)$ be a groupoid object in the category of algebraic spaces
with big fppf stack quotient $[Z/R]$.  Write $\sigma,\tau:R\to Z$ for
the two structural morphisms.  At the moment, we make no (e.g.,
flatness or finiteness) assumptions about $\sigma$ and $\tau$.  
One way to understand the stack $[Z/R]$
is as the stackification (see Lemma 3.2 of \cite{Lau-MB00}) of an
intermediate prestack, which we will denote $\{Z/R\}$.  The objects of
$\{Z/R\}$ over $T$ are given by the elements of $Z(T)$.  Given two such
objects $a,b\in Z(T)$, we define the sheaf of isomorphisms
$\isom_T(a,b)$ to be the fiber of $R\to Z\times Z$ over $(a,b)$. 
Using the groupoid structure on $(R,Z)$, one can check that this
defines a prestack, and that the natural map to the $2$-categorical
fiber product $R\to Z\times_{\{Z/R\}}Z$ is an isomorphism (see Paragraph
2.4.3 of \cite{Lau-MB00}). For any stack
$\ms Y$, the universal property of stackification says that the
restriction functor $$\hom([Z/R],\ms Y)\to\hom(\{Z/R\},\ms Y)$$ is an
equivalence of groupoids.

Let $R^{(2)}$ denote the fiber product
$R\times_ZR$.  The groupoid structure yields
three morphisms $R^{(2)}\to R$: the two projections $\pr_1$ and
$\pr_2$, and the composition map $m$.  Given a prestack $\ms Y$, an
object $\phi:Z\to\ms Y$, and an isomorphism
$\eta:\phi\sigma\simto\phi\tau$, the \emph{coboundary\/} of $\eta$ is
defined to be the element
$$\partial\eta:=(\eta\pr_1)(\eta m)^{-1}(\eta\pr_2)$$ of
$\aut(\phi\sigma\pr_2)$.  

\begin{defn} 
  Let $\ms Y$ be a prestack. Then an \emph{$R$-equivariant object of $\ms
    Y$ over $Z$\/} is an object $\phi:Z\to\ms Y$ and an isomorphism
  $\eta:\phi \sigma\simto\phi \tau$ of morphisms $R\to\ms Y$ such that
  the coboundary $\partial\eta$ is trivial.
\end{defn}

The $R$-equivariant objects of $\ms Y$ naturally form a groupoid, which we will
denote by $\ms Y_{(R,Z)}$.  A $1$-morphism of groupoids $\ms Y\to\ms Y'$ induces a
functor $\ms Y_{(R,Z)}\to\ms Y'_{(R,Z)}$.

A basic example of an equivariant object of a prestack comes from the morphism
$\phi:Z\to \{Z/R\}$ induced by the point $\id\in Z(Z)$.  The isomorphism
$\eta:\phi\sigma\simto\phi\tau$ arises as follows: by definition, we have that
$\isom_R(\phi\sigma,\phi\tau)$ is the fiber product sheaf
$$
\xymatrix{ 
  \isom_R(\phi\sigma,\phi\tau)\ar[r]\ar[d] & R\ar[d]\\
  R\ar[r] & Z\times Z,
}
$$ where both maps $R\to Z\times Z$ are the pair
$(\sigma,\tau)$.  The diagonal of $R\times R$ yields a canonical section of
$\isom_R(\phi\sigma,\phi\tau)$, giving rise to an equivariant object.

In fact, this is the universal equivariant object, as we now make precise.  Given a
prestack $\ms Y$, the constructions of the two previous paragraphs yield a functor
$$\eps:\hom(\{Z/R\},\ms Y)\to\ms Y_{(R,Z)}$$ between groupoids.

\begin{prop}\label{P:equivariant objects} The functor $\eps$ is
  an equivalence for any stack $\ms Y$.
\end{prop}
\begin{proof} 
  We first describe the groupoid $\hom(\{Z/R\},\ms Y)$. Let $\ms P$ be
  the groupoid of pairs $(\phi,\iota)$ consisting of a $1$-morphism
  $\phi:Z\to\ms Y$ and a morphism $\iota:(R,Z)\to (Z\times_{\ms
    Y}Z,Z)$ of groupoids.  The isomorphisms in $\ms P$ are given by
  isomorphisms between the maps $\phi$ which are compatible with the
  maps $\iota$.

Given a morphism
  $\phi:\{Z/R\}\to\ms Y$, composition with the natural morphism $Z\to \{Z/R\}$
  defined above yields a diagram
  $$
  \xymatrix{\hskip-1em R= 
    Z\times_{\{Z/R\}} Z\ar@<-.5ex>[d]\ar@<.5ex>[d]\ar[r] & Z\times_{\ms Y} Z\ar@<-.5ex>[d]\ar@<.5ex>[d]\\
    Z\ar[d]\ar^{\id_Z}[r] & Z\ar[d]\\
    \{Z/R\}\ar[r] & \ms Y.}
  $$ 
  The diagram induces a morphism of groupoids $\iota:(R,Z)\to
  (Z\times_{\ms Y}Z,Z)$, yielding a functor from
  $\hom(\{X/R\},\ms Y)$ to $\ms P$.

  We can produce a functor $\ms P\to\hom(\{Z/R\},\ms Y)$ in the
  opposite direction as follows.  Given a $1$-morphism $\phi:Z\to\ms Y$ and
  a morphism $\iota:(R,Z)\to (Z\times_{\ms Y}Z,Z)$ of groupoids, we
  make a $1$-morphism of prestacks $\{Z/R\}\to\ms Y$ as follows.  An
  object $\alpha\in\{Z/R\}(T)=Z(T)$ gets sent to $\phi\alpha\in\ms
  Y_T$, and an arrow $r\in\isom_{\{Z/R\}(T)}(\alpha,\beta)=R(T)$ gets
  sent to the arrow $\psi:\phi\alpha\simto\phi\beta$ determined by the
  image of $r$ in $Z\times_{\ms Y}Z$ (i.e., $r$ maps to the triple
  $(\alpha,\beta,\psi)$ in the functorial construction of the
  $2$-fiber product).

  The result is an equivalence of groupoids $\hom(\{X/R\},\ms Y)\to\ms P$.
  Sending a pair $(\phi,\iota)$ to the pair
  $(\phi,\iota(\sigma,\tau))$ gives a functor $e:\ms P\to\ms
  Y_{(R,Z)}$ which factorizes $\eps$. The result is thus proven if we
  show that $e$ is an equivalence of groupoids.  A morphism $(R,Z)\to
  (Z\times_{\ms Y}Z,Z)$ extending $\phi$ is given by a morphism
  $y:R\to Z\times_{\ms Y}Z$ with image $(\sigma,\tau)$ in $Z(R)\times
  Z(R)$ such that the composition arrow $R^{(2)}\to R$ is compatible
  via $y$ with the canonical morphism
  $$
  (Z\times_{\ms Y}Z)\times_{Z}(Z\times_{\ms Y}Z)\simto Z\times_{\ms Y}(Z\times_Z
  Z)\times_{\ms Y}Z\simto Z\times_{\ms Y}Z.
  $$ 
  The arrow $y:R\to Z\times_{\ms Y}Z$ gives a triple $(\sigma,\tau,\eta)$ with
  $\eta:\phi\sigma\simto\phi\tau$.  The coboundary condition on $\eta$ is precisely
  the condition that a triple $(\sigma,\tau,\eta)$ give rise to a morphism of
  groupoids over $\phi$, as desired.
\end{proof}

\begin{remark}
  A similar result is proven in Section 3.8 of \cite{MR2223406}, where
  equivariant objects of stacks are treated.  There, the groupoid
  $(R,Z)$ is given by a group action $G\times Z\to Z$.
\end{remark}

In particular, we may apply \eqref{P:equivariant objects} to the case of a group $G$
acting on an algebraic space $Z$, yielding an equivalence between $G$-equivariant
maps $Z\to\ms Y$ and morphisms $[Z/G]\to\ms Y$.  We can use this to prove a purity
theorem for maps $[Z/G]\to\ms Y$.

\begin{prop}\label{P:equivariant purity} Suppose $\ms M$ is a 
  separated Deligne-Mumford stack with coarse moduli space $M$.  Let $(R,Z)$ be a
  groupoid of algebraic spaces with $Z$ regular and $R$ normal and with flat
  structural morphisms.  Suppose $\psi:Z\to M$ is an $R$-invariant morphism,
  $U\subset Z$ is a dense $R$-invariant open subspace, and $\phi:U\to\ms M$ is an
  $R$-equivariant object covering $\psi|_U$.  If $Z\setminus U$ has codimension at
  least $2$ in $Z$ then $\phi$ extends to an $R$-equivariant object of $\ms M$ over
  all of $Z$ which covers $\psi$.
\end{prop}
\begin{proof} By the Purity Lemma \cite[2.4.1 and 2.4.2]{MR1862797},
  $\phi$ lifts to $\widebar\phi:Z\to\ms M$.  It remains to show
  equivariance.  We are given an isomorphism $\alpha:\phi
  \sigma\simto\phi \tau$.  As $\ms M$ is separated,
  $\isom_R(\widebar\phi \sigma,\widebar\phi \tau)$ is finite over $R$
  (via either projection). Furthermore, since $R$ is normal, any
  finite birational morphism $Y\to R$ is an isomorphism.  It follows
  that taking the closure of $\alpha$ in $\isom_R(\widebar\phi
  \sigma,\widebar\phi \tau)$ yields a global section over $R$.
  Moreover, we know that the coboundary of $\alpha$ is trivial over
  the preimage of $U$.  
Since $U$ is schematically dense in $Z$ and
  the structural morphisms of the groupoid are flat, it follows that the
  preimage of $U$ in $R^{(2)}$ is schematically dense.  (To prove
  this, first note that it suffices to prove that some open subspace
  of $U$ has schematically dense preimage in $R^{(2)}$.  Since $Z$ is
  regular and quasi-separated, by working with one component at a time
  we can choose a dense open subspace $U'$ of $U$ whose inclusion
  $i:U'\inj Z$ is a quasi-compact morphism, so that $\ms O_Z\to
  i_\ast\ms O_{U'}$ is an injective map of quasi-coherent $\ms
  O_Z$-algebras.  Since pushforward and flat base change commute for
  quasi-coherent sheaves, we see that the induced map $\ms
  O_{R^{(2)}}\to (i\times\id_{R^{(2)}})\ast\ms O_{U'\times_Z
    R^{(2)}}$ is injective.  This shows that $U'\times_ZR^{(2)}$ is
  schematically dense in $R^{(2)}$, as desired.) 
Using the fact
  that $\isom_R(\widebar\phi\sigma,\widebar\phi\tau)$ is separated
  over $R$, we see that the coboundary of $\alpha$ is trivial over all
  of $R^{(2)}$.
\end{proof}

\begin{cor}\label{C:fluffy lift} 
  Let $\ms Z$ be a smooth Deligne-Mumford stack and $\ms U\subset\ms
  Z$ an open substack of complimentary codimension at least $2$.  Let
  $\ms M$ be a separated Deligne-Mumford stack with coarse moduli
  space $M$.  Given a morphism $\psi:\ms Z\to M$ and a lift $\phi_{\ms
    U}:\ms U\to\ms M$, there is a unique extension $\phi:\ms Z\to\ms
  M$ up to unique isomorphism.
\end{cor}
\begin{proof} We include an alternative, purely stack-theoretic proof
  (without invoking group\-oids).  This proof has the advantage of
  greater intrinsic clarity, although we find the groupoid formalism
  helpful in the proof of \eqref{P:extension} below.

  Consider the morphism $\phi:\ms Z\times_{M}\ms M\to\ms Z$.  By
  assumption, there is a section $\sigma$ over $\ms U$.  Let $\ms
  Y=\widebar{\sigma(\ms U)}$ be the stack-theoretic closure.  The
  projection $\ms Y\to\ms Z$ is proper, quasi-finite, and an
  isomorphism in codimension $1$.  This persists after any \'etale
  base change $Z\to \ms Z$, whence, since $\ms Z$ is smooth, we see
  that $\rho:\ms Y':=\ms Y\times_{\ms Z} Z\to Z$ must be an
  isomorphism.  Indeed, it immediately follows that (via $\rho$) $Z$ is the
  coarse moduli space of $\ms Y'$. 
  On the other hand, over the strict
  localizations of $Z$, $\ms Y'$ is a finite group quotient $[\spec
  R/G]$ with coarse space $\spec S$.  By assumption, $S$ is regular
  and $S\subset R$ is finite and unramified in codimension $1$, hence
  is finite \'etale by purity.  It follows that $\spec S\iso[\spec
  R/G]$.  We conclude that $\ms Y\to\ms Z$ is an isomorphism,
  and thus that there is a lift $\ms Z\to\ms M$ over $M$.
\end{proof}

\subsection{Stacky branched covers}
We briefly recall the basic facts concerning stacky branched covers.  Let $D\subset
Z$ be an effective Cartier divisor in an algebraic space, corresponding to a pair
$(\sL,s)$ with $\sL$ an invertible sheaf on $Z$ and $s\in\Gamma(Z,\sL)$ a regular
global section (i.e., $s$ is not a zero divisor).
Let $\LL$ be the Artin stack $[\A^1/\G_m]$; $\LL$ represents the stack (on the
category of algebraic spaces) of pairs $(\sL,s)$ consisting of an invertible sheaf
and a global (not necessarily regular) section.  The map $x\mapsto x^n$ defines a
morphism $\nu_n:\LL\to\LL$.

\begin{prop}\label{P:stack of roots}  Let $Z,D,\sL,s$ be as above.  
  Define $Z[D^{1/n}]$ to be $Z\times_{(\sL,s),\LL,\nu_n}\LL$.
  \begin{enumerate}
  \item $\pi:Z[D^{1/n}]\to Z$ is a tame Artin stack with coarse moduli
    space $Z$; the natural morphism $Z[D^{1/n}]\times_Z(Z\setminus
    D)\to Z\setminus D$ is an isomorphism.\label{roots1}

  \item $(Z[D^{1/n}]\times_ZD)_{\red}\to D_{\red}$ is the $\m_n$-gerbe
    of $n^{\text{th}}$ roots of the invertible sheaf $\sL|_D$.\label{roots2}

  \item There exists a pair $(\mcL,\sigma)$ of an invertible sheaf and a global
    section on $Z[D^{1/n}]$ with an isomorphism $\mcL^{\tensor
      n}\simto\pi^{\ast}\sL$ sending $\sigma^{\tensor n}$ to $\pi^{\ast}s$. The
    section $\sigma$ is regular.  Moreover, the pair $(\mcL,\sigma)$ is universal:
    $Z[D^{1/n}]$ represents the stack of such pairs of $n^{\text{th}}$
    roots.\label{roots3}

  \item Zariski locally on $Z$, $Z[D^{1/n}]$ has the form $[\rspec(\ms
    O_Z[z]/(z^n-t))/\m_n]$, where $t=0$ is a local equation for
    $D$.\label{roots4}

  \item The stack $Z[D^{1/n}]$ is a global quotient of the form $[Q/\G_m]$. If $Z$
    and $D$ are regular, then so is $Q$.\label{roots5}

  \item If $Z$ is projective over a field and $D$ is smooth then
    there is a finite flat morphism $Y\to Z[D^{1/n}]$ with $Y$ a
    projective scheme.\label{roots6}
  \end{enumerate}
\end{prop}
\begin{proof} The proof of (\ref{P:stack of roots}.\ref{roots1}) through
  (\ref{P:stack of roots}.\ref{roots4}) has been treated numerous times in the
  literature (see for example \cite[4.1]{matsuki-olsson} and \cite[\S 2]{MR2306040}).
  The penultimate statement may be proven as follows: given the universal pair
  $(\mcL,\sigma)$, let $Q\to Z[D^{1/n}]$ be the total space of the $\G_m$-torsor
  associated to $\mcL$.  Since the stabilizer action on each geometric fiber of
  $\mcL$ is faithful, it is a standard result that $Q$ (which is the bundle of frames
  of the line bundle associated to $\mcL$) is an algebraic space.  It immediately
  follows that $Z[D^{1/n}]\iso [Q/\G_m]$.

  To prove the final statement, we recall Viehweg's formulation of the
  Kawamata covering trick \cite[Lemma 2.5]{Viehweg95} and point out a
  slight modification.  Write $Z^{\text{\rm sm}}$ for the smooth locus
  of $Z$; this is an open subscheme containing $D$.  Let $d=\dim Z$.
  Let $H$ be an ample divisor on $Z$.  For sufficiently large $m$, the
  divisor $nmH-D$ is very ample.  Choose general members
  $E_1,\ldots,E_d$ such that $(E_1+E_2+\cdots+E_d+D)|_{Z^{\text{\rm
        sm}}}$ is a \sncd.  Each $E_i+D$ is in $n\Pic$, so we
  can construct the usual cyclic cover branched over $E_i+D$ (see for
  example \cite[Definition 2.49(3)]{MR1658959}), say $C_i\to Z$.  By
  construction, $C_i\to Z$ is a finite flat morphism.  Let $\widebar
  Y:=C_1\times_ZC_2\times_Z\cdots\times_Z C_d$.  Over
  $Z^{\text{\rm sm}}\setminus D$, the transversality of
  $E_1,\ldots,E_d$ ensures that $\widebar Y$ is smooth.  On the
  other hand, one can check that the normalization $Y^{\nu}$ of
  $\widebar Y|_{Z^{\text{\rm sm}}}$ is smooth.  Moreover, the reduced
    structure on the preimage of $D$ in $Y^{\nu}$ gives an effective
    Cartier divisor $D'$ such that $nD'=D|_{Y^{\nu}}$.  Gluing
    $\widebar Y|_{Z\setminus D}$ to $Y^{\nu}$ yields a finite flat
    morphism $f:Y\to Z$ such that $f^{-1}(Z^{\text{\rm sm}})$ is
    smooth and there is an effective Cartier divisor $D'\in Y$ such
    that $nD'=f^{\ast}D$.  By the universal property of $Z[D^{1/n}]$,
  there is a $Z$-morphism $f':Y\to Z[D^{1/n}]$.  Over the complement
  of $D$, $f'$ and $f$ are naturally isomorphic.  On the other hand,
  in a neighborhood of $D$, both $Z[D^{1/n}]$ and $Y$ are regular and
  equidimensional of the same dimension.  Applying \cite[Corollary to
  Theorem 32.1]{MR1011461} to the pullback of $f'$ over an affine
  \'etale neighborhood of $D$, we see that $f'$ is finite and flat, as
  desired.
\end{proof}

Given a \sncd $D=D_1+\cdots+D_{\ell}$ in $Z$ (which implies that the
strict local rings of $Z$ are regular at each point in the support of
$D$), define
$$
Z\langle D^{1/n}\rangle:=Z[D_1^{1/n},\ldots,D_{\ell}^{1/n}]:=
Z[D_1^{1/n}]\times_Z\cdots\times_Z Z[D_{\ell}^{1/n}].
$$ 
In Cadman's notation~\cite{MR2306040}, $X[D^{1/n}]$ is written as $X_{\sL,s,n}$ and
$X\langle D^{1/n}\rangle$ is written as $X_{(D_1,\ldots,D_n),(n,\ldots,n)}$.

We assume in what follows that $Z$ is excellent.  (By definition, an
algebraic space is excellent if \emph{every\/} \'etale cover by a scheme is
excellent. Simply requiring it for one cover is not sufficient, as
shown in \cite[18.7.7]{EGAIV4}.)
\begin{lem} 
  $Z\langle D^{1/n}\rangle$ is regular in a neighborhood of $Z\langle
  D^{1/n}\rangle\times_Z D$.
\end{lem}
\begin{proof} 
  It is clear that the formation of $Z\langle D^{1/n}\rangle$ and the
  statement of the lemma are compatible with \'etale base change, so
  we may assume that $Z=\spec R$ is an affine scheme.  Since $D$ is a
  \sncd and $Z$ is excellent, $Z$ is regular in a Zariski neighborhood
  of $D$.  Upon replacing $Z$ by this neighborhood, we may assume that
  $Z$ is regular.  Shrinking $Z$ further if necessary, we may also
  assume that $\ms O(D_1),\ldots,\ms O(D_{\ell})$ are trivial.  Let
  $t_i=0$ be an equation for $D_i$.  In this case, $Z\langle
  D^{1/n}\rangle$ is locally isomorphic to the quotient stack for the
  action of $\m_n^{\oplus\ell}$ on $Y=\spec
  R[z_1,\ldots,z_{\ell}]/[(z_1^n-t_1,\ldots,z_{\ell}^n-t_{\ell})]$.
  Since $D$ is a \sncd, it is easy to see that $Y$ is regular (and
  excellent).
\end{proof}

\begin{lem} Let $D\subset Z\to S$ be a flat relative \sncd.  The
  formation of $Z\langle D^{1/n}\rangle$ is compatible with base
  change.
\end{lem}
\begin{proof} This follows immediately from the definition.
\end{proof}

\subsection{Applications to lifting problems}

In this section we fix a discrete valuation ring $R$ with uniformizer
$t$, fraction field $K$, and residue field $\rf$.  Let $\ms Z$ be a
tame separated Deligne-Mumford stack of finite type over $R$ with
coarse moduli space $\spec R$ and trivial generic stabilizer.

\begin{lem}\label{L:local splitting} With the above notation, there exists
  a positive integer $n_0$ such that for all $n$ divisible by $n_0$,
  there is a unique $R$-morphism $\spec R[t^{1/n}]\to\ms Z$ up to
  unique isomorphism.
\end{lem}
\begin{proof} Let $R'$ be the strict Henselization of $R$.  By the
  local structure theory of Deligne-Mumford stacks, there exists a
  finite generically Galois extension $S/R'$ with Galois group $G$ and
  an $R'$-isomorphism $\ms Z\tensor R'\iso[S/G]$.  Since $\ms Z$ is
  tame, we may assume that the order of $G$ is invertible in $R$.  Let
  $\widetilde S$ denote the normalization of $S$.  Since $\widetilde
  S/R'$ is generically Galois of degree invertible in $R$, it follows
  from Abhyankar's Lemma and the structure of finite unramified
  extensions of Henselian local rings that $\widetilde S$ is a finite
  product of rings of the form $R'[t^{1/n_0}]$ for a fixed $n_0$. (A
  treatment of the case of discrete valuation rings, which is all we
  use here, may be found in \cite[Chapter IV, \S\S{1-2}]{MR554237}.)
  It follows that the stabilizer of the factor $R'[t^{1/n_0}]$ is
  isomorphic to $\m_{n_0}(\bar\rf)$.  Let $n$ be any integer divisible
  by $n_0$. Then there exists a natural $\m_n$-equivariant map $\spec
  R[t^{1/n}]\to\spec R[t^{1/n_0}]$.  Since $R'[t^{1/n}]=R'\tensor_R
  R[t^{1/n}]$, it follows that there is an \'etale surjection
  $U\to\spec R[t^{1/n}]$ and an $R$-morphism $\phi:U\to\ms Z$.  Let
  $Y=\spec R[t^{1/n}]$.  Over the generic fiber of $U\times_Y U$ there
  is a descent datum $\psi:\pr_1^{\ast}\phi\simto\pr_2^{\ast}\phi$
  arising from the fact that $\ms Z$ is generically isomorphic to $R$.
  Since $U\to Y$ is unramified, $U\times_Y U$ is a Dedekind scheme,
  and since $\ms Z$ is separated, it follows that $\psi$ extends to a
  descent datum for the covering $U\to Y$, yielding an $R$-morphism
  $Y\to\ms Z$.  Uniqueness follows from separatedness and the fact
  that $\ms Z\to\spec R$ is a generic isomorphism.
\end{proof}

\begin{lem}\label{L:codim 1 structure} In the situation of the
  previous lemma, let $Y=\spec R[t^{1/n_0}]$.  There is an induced
  morphism of $R$-stacks $[Y/\m_{n_0}]\to\ms Z$ which identifies
  $[Y/\m_{n_0}]$ with the normalization of $\ms Z$.
\end{lem}
\begin{proof} 
  The scheme $Y\times\m_{n_0}$ is Dedekind, and the generic morphism $Y\tensor
  K\to\ms Z\tensor K$ is clearly equivariant.  Arguing as in the proof of
  \eqref{L:local splitting}, the induced generic isomorphism between the two maps
  $Y\times\m_{n_0}\to\ms Z$ extends (uniquely) over all of $Y\times\m_{n_0}$.  It
  follows that $Y\to\ms Z$ is equivariant, yielding a morphism
  $\rho:[Y/\m_{n_0}]\to\ms Z$ by \eqref{P:equivariant objects}.  Since, in the
  notation of \eqref{L:local splitting}, $\m_{n_0}(\bar\rf)$ is a subgroup of $G$
  (namely, the subgroup fixing the closed point of $\widetilde S$), it follows that
  the morphism of stabilizers induced by $[Y/\m_{n_0}]\to \ms Z$ is injective.  Thus,
  $[Y/\m_{n_0}]\to\ms Z$ is proper, quasi-finite, birational, and injective on
  geometric stabilizers, which implies that it is a finite (affine) morphism.  The
  result follows from the uniqueness of normalization.
\end{proof}

\begin{remark} 
  It is an amusing exercise to understand how \eqref{L:codim 1
    structure} applies to the case of the stack $\ms Z$ given by the
  quotient of a wedge of $n$ lines (in the sense of topology) by one
  of the natural actions of $\Z/n\Z$.  It is easy to see that the
  coarse moduli space is a line, and that there is a single stacky
  point (corresponding to the point at which all of the lines are
  wedged, which is fixed by $\Z/n\Z$).  In particular, $\ms Z$ is
  integral (but admits a finite \'etale cover by a connected reducible
  scheme).  In this case $n_0=1$ and the normalization is simply a
  (non-stacky) line.
\end{remark}

\begin{prop}\label{P:extension} Let $\ms M$ be a proper tame 
  Deligne-Mumford stack with coarse moduli space $M$.  Let $U\subset Z$ be the
  complement of a \sncd $D$ in an excellent scheme.  Suppose $\phi:U\to\ms M$ is a
  morphism such that $\pi\circ\phi:U\to M$ extends to a morphism $\psi:Z\to M$.  Then
  there is an extension of $\phi$ to a morphism $\widetilde{\psi}:Z\langle
  D^{1/N}\rangle\to\ms M$ lifting $\psi$, where $N$ is the least common multiple of
  the orders of the geometric stabilizers of $\ms M$.
\end{prop}
\begin{proof} Let the generic points of $Z\setminus U$ be
  $p_1,\ldots,p_r$.  We claim that it is sufficient to extend $\phi$
  across the Zariski localizations $Z_{p_i}\langle D_i^{1/N}\rangle$
  at each $p_i$.  Indeed, let $Y\to Z\langle D^{1/N}\rangle$ be the
  $\G_m^{\oplus r}$-torsor defined in (\ref{P:stack of
    roots}.\ref{roots5}).  Giving an extension of $\phi$ to $Z\langle
  D^{1/N}\rangle$ is the same as giving an equivariant extension of
  $\phi|_{U_Y}$ to a morphism $Y\to\ms M$.  Such an extension is
  unique up to unique isomorphism, so it immediately follows that once
  one has an extension over each $p_i$, one gets an extension over an
  open subspace $V\subset Z$ whose complement has codimension at least
  $2$ and is contained entirely in the regular locus of $Z$.  Applying
  \eqref{P:equivariant purity} yields the result.

  Thus, let $R$ be the local ring at some $p_i$.  Consider $\ms
  Y:=\spec R\times_M\ms M\to\spec R$.  Since $R$ is normal, we may
  apply \eqref{L:codim 1 structure} to conclude that $\ms Y=\spec
  R[p_i^{1/m}]$, where $m$ is the order of the geometric stabilizer
  over $p_i$.  Since $\spec R[p_i^{1/N}]$ naturally dominates $\spec
  R[p_i^{1/m}]$ over $R$ for any multiple $N$ of $m$, we see that we
  can extend $\phi$ across the preimage of $p_i$ in $Z\langle
  D^{1/N}\rangle$, as required.
\end{proof}

\begin{remark}
  The reason we call this section ``weak stacky stable reduction'' is the following:
  given a discrete valuation ring $R$ and a family over its generic point, the
  methods of this section produce a family over a stack with coarse moduli space
  $\spec R$, \emph{as long as we already have the extension of the coarse moduli
    map\/}.  This makes the statement easier to prove but far weaker than stable
  reduction, even in a stacky form (cf.\ \cite{MR2051053}).
\end{remark}

\section{Proof of the main theorems}

\subsection{Terminology}\label{SS:main}
Let $\ms M^{\circ}$ be a separated Deligne-Mumford stack of finite
type over $T$ with coarse moduli space $\mc M^\circ$.

\begin{defn}
  The stack $\ms M^\circ$ is \emph{compactifiable\/} if there is an
  open immersion $\ms M^\circ\inj\ms M$ into a proper Deligne-Mumford
  stack.  If $\ms M^\circ$ is provided with a compactification, we
  will say it is \emph{compactified\/}.
\end{defn}

\begin{lem}\label{L:quotient-compac}
  Any separated Deligne-Mumford stack arising as the quotient of an
  action by a linear group on a quasi-projective scheme over a field
  is compactifiable.
\end{lem}
\begin{proof}
  By Theorem 5.3 of \cite{kresch-geometry-of-stacks}, such a stack
  admits a locally closed immersion into a smooth proper
  Deligne-Mumford stack with projective coarse moduli space.  Taking
  the stack-theoretic closure yields the result.
\end{proof}

\begin{defn}
  A \emph{coarse compactification\/} of $\ms M^\circ$ is a
  compactification of the coarse moduli space $\mc M^\circ$.  If $\ms
  M^\circ$ is provided with a coarse compactification, we will say it
  is \emph{coarsely compactified\/}.
\end{defn}

Let $\ms M^\circ$ be a coarsely compactified separated Deligne-Mumford
stack of finite type over $T$.  Suppose the coarse compactification
$\mc M^\circ\inj\mc M$ is relatively polarized by $\ms L$.

\begin{defn}\label{D:weakly-bdd-wrt}
  Given a function $\wb:\Z_{\geq 0}^2\to\Z$, we will say that $\ms M^\circ$ is
  \emph{weakly bounded with respect to $\mc M$ and $\ms L$ by $\wb$\/} if for every
  geometric point $t\to T$ and every $(g,d)$-curve $C^\circ\subseteq C$ over
  $\kappa(t)$, every morphism $\xi:C^\circ\to \mc M_t$ factoring through $\ms
  M_t^\circ$ satisfies $\deg \xi_C^*\ms L\leq \wb(g,d)$, where $\xi_C$ is the
  extension of $\xi$ to a morphism $C\to \mc M_t$.
  Cf.\ \eqref{D:coarse weak boundedness}.
\end{defn}

\begin{defn}\label{D:weakly-bdd}
  The stack $\ms M^\circ$ is \emph{weakly bounded} if there exists a
  coarse compactification $\mc M^\circ\inj\mc M$, a relative
  polarization $\ms L$ of $\mc M$ with respect to $\mc M^\circ$, and a
  function $\wb:\Z_{\geq 0}^2\to\Z$ such that $\ms M^\circ$ is weakly
  bounded with respect to $\mc M$ and $\ms L$ by $\wb$.
\end{defn}

Given a scheme $U$, define a relation on the set of isomorphism
classes of morphisms $\phi:U\to\ms M^{\circ}$ as follows:
$\phi_1\sim\phi_2$ if and only if there exists a connected $k$-scheme
$T$, two points $t_1,t_2\in T(k)$, and a morphism $\psi:U\times
T\to\ms M^{\circ}$ such that $\psi|_{U\times t_i}\iso\phi_i$ for
$i=1,2$.  This generates an equivalence relation $\equiv$.

\begin{defn} The equivalence classes for the equivalence relation
  $\equiv$ are called \emph{deformation types\/}.
\end{defn}

It is clear that this notion agrees with \eqref{D:defo types} when $\ms M^{\circ}$ is
the moduli stack of canonically polarized manifolds.

\subsection{The main theorem}\label{SS:mainproof}

\begin{proof}[Proof of Theorem \ref{T:main theorem}]
  Observe that by \eqref{C:cc-quasi-proj} there is a finite type extension
  $\widetilde T\to T$ and a morphism $U_{\widetilde T}\to\mc M^\circ$ with the
  following property:  For a geometric point $t\to T$, every morphism
  $U_t\to\mc M^\circ_t$ that arises by composition $U_t\to\ms
  M^\circ\to\mc M^\circ$ is parametrized by a point of $\widetilde T$.
  Let $\ms M^\circ\inj\ms M$ be a compactification of $\ms M^\circ$,
  and let ${\mc M}$ be the coarse moduli space of $\ms M$.  Applying
  \eqref{L:compactifying family}, there is a further finite type
  extension $\sigma:T'\to\widetilde T$, a proper morphism $\mc B\to T'$, a
  relative \sncd $D\subset\mc B$, an isomorphism $\mc B\setminus D\iso
  U_{T'}$, and a morphism $\mc B\to{\mc M}$ such that for every fiber
  $U_t$ and every morphism $\phi:U_t\to\ms M^{\circ}$, there exists a
  point $t'\to T'_t$ such that $\sigma(t')=t$ and the restriction of the induced morphism
  $\mc B_{t'}\to{\mc M}_t$ to $U_{t'}$ 
  is the coarse morphism associated to $\phi$.

  By \eqref{P:extension} and the fact that $T$ has characteristic $0$, the morphism
  $U_{t'}\to\ms M^{\circ}_{t'}$ extends to a morphism $\mc B\langle
  D^{1/N}\rangle_{t'}\to\ms M_{t'}$ over the coarse moduli map $\mc B_{t'}\to{\mc
    M}_{t'}$ for any geometric point $t'\to T'$.  Consider the morphism of stacks
  $$\mu:\uhom_{T'}(\mc B\langle D^{1/N}\rangle,\ms M_{T'})\to\uhom_{T'}(\mc B,{\mc
    M}_{T'}).$$ We know that $\mu$ is of finite type: if $U$ is
  quasi-projective, then this follows from \cite[Theorem
  C.4]{abramovich-olsson-vistoli}, as all of the stacks involved are
  tame (the characteristic of $k$ being $0$) and separated, and $\mc
  B\langle D^{1/N}\rangle$ is proper and flat over the base.
  Furthermore, by \eqref{L:inclusion} there is a finite type
  monomorphism
$$\ms S\to\uhom_{T'}(\mc B\langle
D^{1/N}\rangle,\ms M_{T'})$$ parametrizing morphisms that pull back
the boundary $\ms M\setminus\ms M^{\circ}$ into $D$.  The given
``universal'' coarse moduli map $\mc B\to {\mc M}$ determines a
section of $\uhom_{T'}(\mc B,{\mc M}_{T'})$ over $T'$.  Pulling this
back to $\ms S$ yields a finite type $T'$-stack $\ms H\to T'$ such
that for any geometric point $t\to T$, the set of deformation types of
morphisms $U_t\to\ms M^{\circ}$ is a quotient of the set of connected
components of $\ms H_t$ (the fiber of the morphism $\ms H\to T'\to
T$).  Indeed, any point of $\ms H_t$ parametrizes a morphism sending
$U_t$ to $\ms M^\circ$ and any such morphism occurs as such a point,
so any two points in the same connected component represent
deformation equivalent morphisms $U_t\to\ms M^\circ$, and any
deformation type is represented by a point of $\ms H_t$.  (There could
conceivably be deformation equivalent morphisms which lie in different
components of $\ms H_t$, as our construction makes frequent use of
stratification.)  Since $\ms H\to \mc B$ is of finite type, the number
of connected components is bounded above for all points $t$, giving a
bound on the number of deformation types.
\end{proof}

\begin{cor}\label{C:uniform curve} If $\ms M^{\circ}$ is
  weakly bounded then there exists a function $\wb_{\ms M}:\Z_{\geq
    0}^2\to\Z$ such that for every smooth curve $C$ of genus $g$ with
  $d$ marked points $p_1,\ldots,p_d$, the number of deformation types
  of morphisms $C\setminus\{p_1,\ldots,p_d\}\to\ms M^{\circ}$ is
  finite and bounded above by $\wb_{\ms M}(g,d)$.
\end{cor}
\begin{proof} 
  Choosing an affine cover of $\ms M_{g,d}$ and pulling back the
  universal curve yields a quasi-compact family containing all
  $d$-pointed smooth curves of genus $g$.  The result thus follows
  immediately from \eqref{T:main theorem}.
\end{proof}

\begin{remark}
  The uniformity result of \eqref{C:uniform curve} was first proven by Caporaso for
  families of curves (i.e., for $\ms M^{\circ}=\ms M_q$) in \cite{Caporaso02}, using
  methods specific to the stack of curves.  In \cite{Heier04b}, Heier refined
  Caporaso's results to produce an effective uniform bound.  It would perhaps be
  interesting to determine what auxilliary data about the stack $\ms M^{\circ}$ are
  necessary to prove an abstract effective form of \eqref{C:uniform curve}.
\end{remark}

\section{Finiteness of infinitesimally rigid families}


Let $\ms M^{\circ}$ be a Deligne-Mumford stack and let $N$ be the least common
multiple of the orders of the stabilizers of geometric points of $\ms
M$.  Suppose $U$ is a $k$-scheme.

\begin{defn}
  A morphism $\chi:U\to\ms M^{\circ}$ is \emph{
    infinitesimally rigid\/} if for every $n\geq 0$, any two
  extensions of $\chi$ to $U\tensor_k k[t]/(t^n)$ are isomorphic.
\end{defn}

Since the diagonal of $\ms M^\circ$ is unramified, there is at most
one isomorphism between two extensions of $\chi$.

\begin{thm}\label{T:strongly non-isotrivial}
  Let $U$ be a smooth variety.  If $\ms M^{\circ}$ is a weakly bounded
  compactifiable Deligne-Mumford stack then the set of isomorphism
  classes of infinitesimally rigid morphisms $U\to\ms M^{\circ}$ is
  finite.  Moreover, the number of isomorphism classes is bounded in a
  manner which is uniform in any finite type family of bases $U$.
\end{thm}

In the standard terminology, this theorem says that ``infinitesimal
rigidity implies rigidity.''  For applications of this result to
families of canonically polarized manifolds, see Section
\ref{sec:apple-canon-polar}.

We start with two lemmata.

\begin{lem}
  \label{L:canonical-cover}
  Let $S$ be a reduced locally Noetherian scheme, $\pi:Z\to S$ and
  $P\to S$ two $S$-schemes of finite type with $P$ separated over $S$.
  Further let $V\subset Z$ be an open subscheme that is dense in every
  fiber $Z_s$. Assume that $V\to S$ has a section and the geometric
  fibers of $\pi$ are reduced.  Then
  \begin{enumerate}
  \item any $S$-morphism $\zeta:Z\to P$ such that the restriction of $\zeta$
    to each geometric fiber $V_s$ is constant factors through a
    section $S\to P$;
  \item if in addition $S$ is of finite type over a field, it is
    sufficient for $\zeta$ to be constant on geometric fibers over closed
    points of $S$.
  \end{enumerate}
\end{lem}
\begin{proof}
  The statement is local on $S$, so we may assume that $S$ is
  Noetherian.  Since $V\to S$ has a section, it is a universal
  effective epimorphism \cite[IV.1.12]{SGA3}.  Since $S$ is reduced
  and Noetherian, it has a dense subscheme consisting of finitely many
  reduced points $t_1,\dots,t_n$ (the generic points of the
  irreducible components of $S$).  Extending the residue field of
  $t_i$ is also a universal effective epimorphism, so if $\zeta$ is
  constant on the geometric fiber of $V$ over $t_i$, it must be
  constant on $V\tensor\kappa(t_i)$ for each $i$.  Write $p$ and $q$
  for the two projections $V\times_S V\to V$.  In the exact diagram
  $$
  \hom(S,P)\to\hom(V,P)\rightrightarrows\hom(V\times_S V,P),
  $$
  we have that the two compositions $\zeta p$ and $\zeta q$
  agree on the fibers over each $t_i$.  Since these fibers are dense
  in $V$ and $P$ is separated, the two maps agree on all of $V\times_S
  V$, whence there is a morphism $\gamma:S\to P$ such that
  $\zeta|_V=\gamma\circ\pi|_V$.  Since $V$ is everywhere dense in $Z$
  and $P$ is separated, it follows that $\zeta$ factors through $S$,
  as required.

  The second statement works precisely the same way, using the fact that for any
  closed set $F$ containing all of the closed fibers of $\pi$ we have $F=Z$.
\end{proof}

\begin{lem}\label{L:stupid lemma}
  Let $R$ be a ring and $y\in R$ a regular element.  Let
  $R[\eps]:=R[x]/(x^n)$.  Let $A$ be a finite $R[\eps]$-algebra such
  that the natural maps $R\to A/\eps A$ and
  $R[\eps][1/y]\to A[1/y]$ are isomorphisms.  Then $R[\eps]\to A$ is
  an isomorphism.
\end{lem}
\begin{proof} 
  In the diagram 
  $$
  \xymatrix{%
    R[\eps] \ar[r] \ar_\iota[d] & A \ar[d]\\
    R[\eps][1/y] \ar_\rho[r]  & A[1/y]\\
  }
  $$
  the natural maps $\iota$ and $\rho$ are injective by the hypotheses
  and hence $R[\eps]\to A$ is injective as well.  
  On the other hand, $R[\eps]\to A$ is surjective modulo the nilpotent
  $\eps$, which implies that $R[\eps]\to A$ is itself surjective.
\end{proof}

\begin{proof}[Proof of \ref{T:strongly non-isotrivial}]

  We already know from Theorem \ref{T:main theorem} that the set of deformation types
  of infinitesimally rigid morphisms $U\to\ms M^\circ$ is finite and of cardinality
  bounded above in a finite type family of bases $U$.  To show finiteness of the set
  of infinitesimally rigid morphisms, it thus suffices to show the following: \emph{
    if $(T,t)$ is a pointed smooth connected curve over $k$ and $\Xi^\circ:U\times
    T\to\ms M^{\circ}$ is a morphism whose restriction to $U_t$ is infinitesimally
    rigid then there is a finite base change $T'\to T$ such that
    $\Xi^\circ|_{T'}\cong\Xi^\circ_t\times \id_{T'}$\/}.  We will refer to this
  statement as $(\dagger)$.  If $(\dagger)$ holds then any two deformation equivalent
  infinitesimally rigid morphisms are in fact isomorphic, as desired.

  To show $(\dagger)$, we first note that by \eqref{L:canonical-cover} with $U=S$,
  $V=Z=U\times T$, and $P=U\times M^{\circ}$, the coarse morphism $U\times T\to M^{\circ}$
  factors through a morphism $\chi:U\to M^{\circ}$.  Indeed, it suffices to show that for each
  closed point $u\in U$, the induced map 
  $\widetilde\Xi^\circ_u:T_u\to M^{\circ}$ is constant.  Since $\Xi^\circ_t$ is
  infinitesimally rigid, for every $n\geq 0$, the map
  $(\widetilde\Xi^\circ_u)|_{\spec\ms O_{T_u,t}/\mf m_t^n}$ factors through the
  natural map $\spec\ms O_{T_u,t}/\mf m_t^n\to\spec\ms O_{T_u,t}/\mf m_t$.  It
  follows that the induced map $\spec\widehat{\ms O}_{T_u,t}\to M^{\circ}$ factors
  through the section $t_u\to T_u$, so that $\widetilde\Xi^\circ_u$ satisfies the
  hypotheses of \eqref{L:canonical-cover}.
  
  Choose a compactification of $U\subset B$ by a \sncd $D$ over which there is an
  extension of $\chi$ to a morphism $B\to M$.  By \eqref{P:extension} (using the fact
  that $T$ is regular), we can extend $\Xi^\circ$ to a morphism $\Xi:B\langle
  D^{1/N}\rangle\times T\to \ms M$, corresponding to a homomorphism
  $T\to\uhom(B\langle D^{1/N}\rangle,\ms M)$.  (Note that $\uhom(B\langle
  D^{1/N}\rangle,\ms M)$ is a separated Deligne-Mumford stack by
  \cite[Theorem~C.4]{abramovich-olsson-vistoli} and the fact that $B\langle
  D^{1/N}\rangle$ is smooth.) We claim that any morphism $B\langle
  D^{1/N}\rangle\to\ms M$ whose restriction to $U$ is infinitesimally rigid is
  infinitesimally rigid.

  Granting this claim, let us demonstrate that $(\dagger)$ follows.  An
  infinitesimally rigid point $\xi:\spec k\to\uhom(B\langle D^{1/N}\rangle,\ms M)$
  has the property that any extension of $\xi$ to $\spec k[\eps]$ is isomorphic to
  $\xi\times\spec k[\eps]$, and there is a unique such isomorphism extending the
  identity over the closed point $\spec k\inj\spec k[\eps]$.  In particular, the
  miniversal deformation of $\xi$ is isomorphic to $\spec k$, showing that $\xi$ is a
  smooth morphism.  It follows that the residual gerbe of $\xi$ is a connected
  component of $\uhom(B\langle D^{1/N}\rangle,\ms M)$.  Applied to the morphism
  $\Xi$, this implies that $\Xi|_{\spec\widehat{\ms O}_{T,t}}$ is isomorphic to
  $\Xi_t\times\id_{\spec\widehat{\ms O}_{T,t}}$.

  Now consider the finite scheme $I:=\Isom_T(\Xi,\Xi_t\times \id_T)\to T$.
  By assumption, $I(\widehat{\ms O}_{T,t})\neq\emptyset$.  Applying Popescu's theorem
  (see, e.g., \cite{MR1647069}), the excellence of $\ms O_{T,t}$ implies that
  $\widehat{\ms O}_{T,t}$ is a filtering colimit of smooth $\ms O_{T,t}$-algebras,
  and since $I$ is locally of finite presentation there is thus a smooth $T$-scheme
  $\widetilde T\to T$ such that $I(\widetilde T)\neq\emptyset$.  Since any smooth
  $T$-scheme has \'etale-local sections around any point, we find a quasi-finite
  generically \'etale morphism $T''\to T$ whose image contains $t$ such that
  $I(T'')\neq\emptyset$.  Letting $T'$ equal the normalization of $T$ in the function
  field of $T''$, we find a finite morphism $T'\to T$ such that there is a generic
  isomorphism between $\Xi|_T'$ and $\Xi_t\times \id_{T'}$.  Since
  $\Isom(\Xi|_{T'},\Xi_t\times \id_{T'})\to T'$ is finite and $T'$ is a Dedekind
  scheme, any generic section extends to a global section by the valuative criterion
  of properness.  Restricting to $U\times T'\subset B\langle
  D^{1/N}\rangle\times_TT'$, we see that $(\dagger)$ is established.

  Thus, it remains to show that if $\xi:B\langle D^{1/N}\rangle\to\ms
  M$ maps $U$ to $\ms M^{\circ}$ and $\xi_U$ is infinitesimally rigid
  then $\xi$ itself is infinitesimally rigid.  Let $\xi_1$ and $\xi_2$
  be two infinitesimal deformations of $\xi$ over
  $k[\eps]:=k[x]/(x^n)$.  Consider the sheaf $I:=\isom_{B\langle
    D^{1/N}\rangle[\eps]}(\xi_1,\xi_2)$ on the \'etale site of
  $B\langle D^{1/N}\rangle$.  Since $\ms M$ is separated, $\pi:I\to
  B\langle D^{1/N}\rangle[\eps]$ is a finite representable morphism of
  stacks.  By the definition of infinitesimal rigidity, we have that
  there is a section $\sigma:U[\eps]\to I$ of $\pi$.  

  Let $J\subset I$ be the stack-theoretic closure of
  $\sigma(U[\eps])$, so that $J\to B\langle D^{1/N}\rangle[\eps]$ is
  finite, representable, and an isomorphism over a dense open
  subscheme.  We claim that $J\to B\langle D^{1/N}\rangle[\eps]$ is an
  isomorphism.  To show this, it suffices to work \'etale-locally on
  $B\langle D^{1/N}\rangle$.  Indeed, since $\sigma$ is a
  quasi-compact morphism we have that $J$ is defined by the
  quasi-coherent kernel of the natural morphism $\ms
  O_I\to\sigma_{\ast}\ms O_{U[\eps]}$ of quasi-coherent sheaves.
  Since the formation of this kernel commutes with \'etale base change
  on $B\langle D^{1/n}\rangle$, we see that the formation of $J$
  commutes with \'etale base change on $B\langle D^{1/n}\rangle$.

  Let $R$ be an \'etale local ring of $B\langle D^{1/N}\rangle$, so that $J$ is
  locally represented by an $R[\eps]$-algebra $A$.  Since $R$ is regular and $\sigma$
  is defined over a dense open substack $U$, it follows that $A$ satisfies the
  conditions of \eqref{L:stupid lemma}.  Indeed, if $R\to A/\eps A$ is not an
  isomorphism, then $A/\eps A$ cannot be irreducible (as $R$ is normal and $R\to
  A/\eps A$ is finite and birational).  But $\spec A$ is irreducible, being the
  scheme-theoretic closure of an irreducible scheme.  On the other hand, if $f$ is
  any element of $R$ vanishing on the complement of $U$, we see that $\spec A[1/f]$
  is contained in the open subscheme $\sigma(U)\subset J$, which is isomorphic to $U$
  by definition; thus, $R[\eps][1/f]\to A[1/f]$ is an isomorphism.  We conclude that
  $J$ is an isomorphism, and thus that $\xi_1\iso\xi_2$ via an isomorphism extending
  the given isomorphism over $U[\eps]$, as required.
\end{proof}

\begin{newnumr}
  It may seem that the condition of infinitesimal rigidity is unnatural, especially
  for families over non-proper base varieties $U$.  For families over curves this is
  true (in fact, infinitesimal rigidity almost never holds for families of
  canonically polarized manifolds over an affine curve).  However, for bases $U$ such
  that the boundary divisor in a compactification $B$ is non-ample, there are many
  examples of infinitesimally rigid families.  This is captured in the following
  proposition (which is far from optimal, but serves to illustrate the point).
\end{newnumr}
\begin{prop}
  Let $B$ be a proper smooth $k$-variety and $D\subset B$ a smooth
  irreducible divisor. Assume that $\xi:B\to\ms M$ is an
  infinitesimally rigid morphism to a Deligne-Mumford $k$-stack such
  that $\xi|_D$ is also infinitesimally rigid and that
  $\xi^{\ast}\Omega^1_{\ms M/k}$ is locally free in a neighborhood of
  $D$.  If $\Gamma(D,\ms O(-D)|_D)\neq 0$ then $\xi|_U$ is
  infinitesimally rigid.
\end{prop}
\begin{proof}
  It is a standard fact \cite[III.2.2.2]{MR0491680} that the
  first-order infinitesimal deformations of $\xi$ form a torsor under 
  $\Hom(\xi^{\ast}\operatorname{L}_{\ms M/k},\ms O_B)$, where
  $\operatorname{L}_{\ms M/k}$ is the cotangent complex of $\ms M$ over
  $k$.  Since $L_{\ms M/k}$ is bounded above at $0$ and $\ms H^0(L_{\ms
    M})\cong\Omega^1_{\ms M/k}$, this space is just
  $\Hom(\xi^{\ast}\Omega^1_{\ms M/k},\ms O_B)=\Gamma(B,\ms H)$, where
  $\ms H=\shom(\xi^{\ast}\Omega^1_{\ms M/k},\ms O_B)$.  We know that
  $\Gamma(B,\ms H)=0=\Gamma(D,\ms H|_D)$ and we wish to conclude that
  $\Gamma(U,\ms H)=0$. %
  Any section of $\ms H|_U$ extends to a section of $\ms H(nD)$ for
  some $n\in\bN$, so it suffices to show that $\Gamma(B,\ms H(nD))=0$
  for all $n\in\bN$.
  
  Since $\ms H|_D$ is locally free, any non-zero section of $\ms
  O(-D)|_D$ is $\ms H|_D$-regular.  Thus, since $\Gamma(D,\ms
  H|_D)=0$, it follows that
  $\Gamma(D,\ms H(nD)|_D)=0$ for all $n\in\bN$.  Consider the
  sequence
  $$
  0\to\ms H((n-1)D)\to\ms H(nD)\to\ms H(nD)|_D\to 0.
  $$
  Since $\Gamma(D,\ms H(nD)|_D)=0$, it follows that $\Gamma(X,\ms
  H((n-1)D))\to \Gamma(X,\ms H(nD))$ is an isomorphism for all
  $n\in\bN$, and since $\Gamma(X,\ms H)=0$, it follows that
  $\Gamma(X,\ms H(nD))=0$ for all $n\in\bN$.
\end{proof}
 
An example of this phenomenon arises by considering families over
$C\times C$, where $C$ is a curve of high genus.  If $D\subset C\times
C$ is the diagonal, it follows from the adjunction formula that $\ms
O(-D)|_D\iso\Omega^1_D$, which is globally generated.  If $X\to C$ and
$Y\to C$ are two infinitesimally rigid families of smooth canonically
polarized varieties (e.g., two non-isotrivial families of smooth
curves) then the fiber product $X\times_C Y$ is also infinitesimally
rigid.  Similarly, by the K\"unneth formula, it is easy to see that
the family $Z:=X\times Y\to C\times C$ is infinitesimally rigid.
Applying the proposition, it follows that the restricted family
$Z|_{C\times C\setminus D}$ is infinitesimally rigid.

\section{Applications to canonically polarized
  varieties}\label{sec:apple-canon-polar} 

Write $\ms M^{\circ}_h$ for the (Deligne-Mumford) stack of canonically
polarized manifolds with Hilbert polynomial $h$ and $M^{\circ}_h$ for
its coarse moduli space.  If $\mf f:\ms X\to\ms M^{\circ}_h$ is the
universal family, then the invertible sheaf
$$
\lambda_m^{(p)}:=\left(\det\mf f_{\ast}\omega_{\ms X/\ms
    M^{\circ}_h}^{\tensor m}\right)^{p}
$$ is the pullback of an ample invertible sheaf on
$M^{\circ}_h$ \cite[Theorem 1.11]{Viehweg95}.

We recall a well-known fact about $\ms M^{\circ}_h$. (A similar statement probably
first appeared in a lecture of M.~Artin \cite[2.8]{Kollar90}.)

\begin{lem}\label{L:can-pol-quot}
  The stack $\ms M^{\circ}_h$ is isomorphic to a separated stack of the form
  $[U/\PGL_r]$ (for $r=h(m)$ with $m$ sufficiently large), where $U$
  is a quasi-projective $k$-scheme.
\end{lem}
\begin{proof}
  By Matsusaka's Big Theorem (\cite[Theorem~2]{MR0263816},
  \cite[Theorem~4.2]{Matsusaka72}), there is a positive integer $m$ such that for any
  canonically polarized manifold $X$ with Hilbert polynomial $h$, the global sections
  of $\omega_X^{\tensor m}$ give a non-degenerate embedding into $\P^{h(m)-1}$.  Let
  $H$ be the Hilbert scheme parametrizing closed subschemes of $\P^{h(m)-1}$ with
  Hilbert polynomial $h$; it is well-known that $H$ is projective.  There is an open
  subscheme $V\subset H$ parametrizing closed subschemes which are smooth and
  geometrically connected.  Let $\mc X\to V$ be the universal family with universal
  embedding $\Upsilon:\mc X\inj\P^{h(m)-1}_V$.  Consider the invertible sheaf $\ms
  L:=\Upsilon^{\ast}\ms O(1)\tensor(\omega_{\mc X/V}^{\tensor m})^{\vee}$.  It
  follows from cohomology and base change that there is a closed subscheme $U\subset
  V$ parametrizing the locus over which $\ms L$ is isomorphic to the trivial
  invertible sheaf (see e.g.\ the proof of \cite[Corollary~6 of
  \S{II.5}]{Mumford70}).  It is easy to see that $U$ is $\PGL_{h(m)}$-invariant, and
  it follows from standard methods that $\ms M^{\circ}_h\iso[U/\PGL_{h(m)}]$.

  That $\ms M^{\circ}_h$ is separated follows easily from the fact that any family of
  canonically polarized manifolds is its own relative canonical model. Indeed, using
  the valuative criterion of separatedness the question reduces to following
  statement: If two families of canonically polarized manifolds are given over the
  same smooth curve such that they agree over an open set, then they agree
  everywhere. However, this follows from the fact that within a fixed a birational
  class the relative canonical model over a fixed base is unique. To see this, let
  $f:X\to C$ be one of the families. Then the relative canonical model is
  $\Proj_C\left(\sum_{m\geq 0} f_*\omega_{X/C}^m\right)\to C$. Since $\omega_{X/C}$
  is relatively ample, this is actually isomorphic to $f$. On the other hand, the
  sheaves $f_*\omega_{X}^m$ are birational invariants and since $C$ is fixed, this
  means that so is $\Proj_C\left(\sum_{m\geq 0} f_*\omega_{X/C}^m\right)\to C$.
\end{proof}

\begin{lem}\label{L:can-pol-man-is-compactifiable}
  The stack $\ms M^{\circ}_h$ is weakly bounded and compactifiable.
\end{lem}
\begin{proof}
  The compactifiability of $\ms M^{\circ}_h$ follows from
  \eqref{L:can-pol-quot} and \eqref{L:quotient-compac}.  Weak boundedness
  is much more subtle.  Given $m>0$, Viehweg \cite[Theorem
  3]{Viehweg06a} produced a projective compactification $M_h$ of
  $(M^\circ_h)_{\text{red}}$ and an invertible sheaf
  $\lambda_{m}^{(p)}\in\Pic(M_h)$, nef and ample with respect to
  $(M^{\circ}_h)_{\text{red}}$, such that for any morphism $\xi:C\to
  M_h$ induced by a semistable family $f:X\to C$, we have that
  $\xi^{\ast}\lambda_{m}^{(p)}=\det(f_{\ast}\omega_{X/C}^{m})^p$.

  We claim that $\ms M^{\circ}_h$ is weakly bounded with respect to
  $M_h$ and $\lambda_{m}^{(p)}$, as in \eqref{D:weakly-bdd-wrt}.  The
  proof is similar to the proof of Corollary 4.1 and Addendum 4.2
  of \cite{Bedulev-Viehweg00}.  Let $C^{\circ}\subset C$ be a
  $(g,d)$-curve and let $f^{\circ}:X^{\circ}\to C^{\circ}$ be a family
  of canonically polarized manifolds with Hilbert polynomial $h$.
  There exists a morphism of smooth projective varieties $f:X\to C$
  including $f^{\circ}:X^{\circ}\to C^{\circ}$ as an ``open
  subdiagram.''  By the semistable reduction theorem \cite[Chapter
  II]{KKMS73}, there is a finite morphism $\gamma:D\to C$ and a
  diagram
  $$
  \xymatrix{
    X \ar_f[d] & X_D\iso X\times_C D \ar[l] \ar [d] & Y \ar[l] \ar^{f'}[dl]\\
    C & \ar^\gamma[l] D }
  $$
  with $Y$ semistable over $D$ and $Y\to X_D$ a resolution of
  singularities.  By \cite[Lemma 3.2]{Viehweg83a}, there is an
  inclusion
  $$f'_{\ast}\omega_{Y/D}^{m}\inj\gamma^{\ast}f_{\ast}\omega_{X/C}^{m}.$$
  Thus,
  $$
  \deg \left(\det(f'_{\ast}\omega_{Y/D}^{m})\right)\leq 
  \left(\deg\gamma\right)\deg\left(\det(f_{\ast}\omega_{X/C}^{m})\right).
  $$
  The composed map $D\to M_h$ comes from a semistable family, so that
  (by the result of Viehweg quoted in the previous paragraph) 
  $$\deg(\gamma\circ\xi)^{\ast}\lambda_{m}^{(p)}=\deg
  \left(\det(f'_{\ast}\omega_{Y/D}^{m})^p\right).$$ It follows that
  $$
  \deg\xi^{\ast}\lambda_{m}^{(p)}\leq
  \deg\left(\det(f_{\ast}\omega_{X/C}^{m})^p\right).
  $$
  By \cite[Theorem 1.4(c)]{Bedulev-Viehweg00}, for $m$ sufficiently
  large and divisible the right-hand side of the last equation is
  bounded above by an explicit polynomial in $g,d,n$ and some
  constants depending upon $m$ (which are fixed once $h$ is fixed).
\end{proof}

Theorem~\ref{T:benjamin} is now an immediate corollary of
Theorem~\ref{T:main theorem}.

\begin{remark}
  At the time of this writing, the finite generation of the canonical
  ring has apparently just been proven \cite{math.AG/0610203}.  It has
  been claimed that under the assumption of the minimal model program
  in dimension $\deg h+1$ (in fact, one seemingly needs only the
  existence of relative canonical models), one knows that there is a
  compactification $\ms M^{\circ}_h\subset\ms M_h$ and an invertible
  sheaf $\ms L$ on $\ms M_h$ such that (1) $\ms L|_{(\ms
    M^{\circ}_h)_{\text{\rm red}}}\iso\lambda_m^{(p)}$ for fixed
  sufficiently large and divisible $m$ and $p$, and (2) $\ms L$ is the
  pullback of an invertible sheaf from the coarse moduli space $M_h$
  of $\ms M_h$.  Using these results would give a more natural proof
  of \ref{L:can-pol-man-is-compactifiable}.  Unfortunately, at the
  present time a proper explanation of this implication is not in the
  literature, so we find it prudent to include an alternative proof.
\end{remark}

\begin{remark}
  Because of the terminology that has been used in studying this
  problem, it behooves us to point out that the powerful results of
  Viehweg and Zuo \cite{Vie-Zuo01}, \cite{Vie-Zuo02},
  \cite{Vie-Zuo03a}, concerning the boundedness problem for families
  of varieties, fall short of addressing the entire question.  In
  particular, without the use of stack-theoretic methods, the
  numerical boundedness results (usually referred to as ``weak
  boundedness'') are not enough in themselves to show constructibility
  of the locus of coarse moduli maps arising from families.  It is
  only by combining the numerical results with a study of lifts of
  coarse maps into stacks that one can prove the concrete boundedness
  results of \eqref{T:benjamin} and \eqref{T:main theorem}.  This fact
  is implicit in the work of Caporaso \cite{Caporaso02}, but rather
  than lifting to the stack, she lifted to a level cover of the stack
  of curves. This allowed her to avoid the use of stack-theoretic
  constructions but limited the argument to handle only families of
  curves.
\end{remark}

\begin{cor}\label{C:uniformity over curves} 
  \protect{\rm (cf.\ \cite{Caporaso02,Heier04b} for families of curves)} There exists
  a function $\wb_h(g,d)$ such that for any $d$-pointed smooth projective curve of
  genus $g$, $(C,p_{1},\ldots,p_{d})$, the number of deformation types of families of
  canonically polarized manifolds $X$ over $C\setminus\{p_{1},\dots,p_{d}\}$ with
  Hilbert polynomial $h$ is bounded above by $\wb_h(g,d)$.
\end{cor}
\begin{proof} This is an application of
  \eqref{C:uniform curve}.
\end{proof}


\def\cprime{$'$} \def\polhk#1{\setbox0=\hbox{#1}{\ooalign{\hidewidth
  \lower1.5ex\hbox{`}\hidewidth\crcr\unhbox0}}} \def\cprime{$'$}
  \def\cprime{$'$} \def\cprime{$'$} \def\cprime{$'$}
  \def\polhk#1{\setbox0=\hbox{#1}{\ooalign{\hidewidth
  \lower1.5ex\hbox{`}\hidewidth\crcr\unhbox0}}} \def\cdprime{$''$}
  \def\cprime{$'$} \def\cprime{$'$}
\providecommand{\bysame}{\leavevmode\hbox to3em{\hrulefill}\thinspace}
\providecommand{\MR}{\relax\ifhmode\unskip\space\fi MR}
\providecommand{\MRhref}[2]{%
  \href{http://www.ams.org/mathscinet-getitem?mr=#1}{#2}
}
\providecommand{\href}[2]{#2}

\end{document}